\documentclass[12pt,a4paper]{article}
\usepackage[utf8]{inputenc}
\usepackage[a4paper,margin=2.3cm]{geometry}

\usepackage{verbatim}
\usepackage{amsmath,amssymb,amsthm,mathtools}
\usepackage{aliascnt}
\usepackage{booktabs,longtable,tabularx,array,multirow}
\usepackage{graphicx,subcaption}
\usepackage{enumitem}
\usepackage{xcolor}
\usepackage{listings}
\usepackage[bookmarks=false]{hyperref}
\usepackage[nameinlink,noabbrev]{cleveref}
\usepackage{fancyhdr}

\hypersetup{
  colorlinks=true,
  linkcolor=blue!45!black,
  citecolor=blue!45!black,
  urlcolor=blue!55!black
}

\newtheorem{theorem}{Theorem}[section]
\newaliascnt{lemma}{theorem}
\newtheorem{lemma}[lemma]{Lemma}
\aliascntresetthe{lemma}
\newaliascnt{proposition}{theorem}
\newtheorem{proposition}[proposition]{Proposition}
\aliascntresetthe{proposition}
\newaliascnt{corollary}{theorem}
\newtheorem{corollary}[corollary]{Corollary}
\aliascntresetthe{corollary}
\theoremstyle{definition}
\newaliascnt{definition}{theorem}
\newtheorem{definition}[definition]{Definition}
\aliascntresetthe{definition}
\newaliascnt{remark}{theorem}
\newtheorem{remark}[remark]{Remark}
\aliascntresetthe{remark}
\newaliascnt{example}{theorem}

\aliascntresetthe{example}

\crefname{theorem}{Theorem}{Theorems}
\crefname{lemma}{Lemma}{Lemmas}
\crefname{proposition}{Proposition}{Propositions}
\crefname{corollary}{Corollary}{Corollaries}
\crefname{definition}{Definition}{Definitions}
\crefname{remark}{Remark}{Remarks}
\crefname{example}{Example}{Examples}
\crefname{equation}{equation}{equations}
\crefname{table}{Table}{Tables}
\crefname{figure}{Figure}{Figures}

\newcommand{\K}{\mathcal K}
\newcommand{\U}{\mathcal U}
\newcommand{\PM}{\mathcal M}
\newcommand{\Neu}{\operatorname{Neu}}

\newcommand{\Act}{\operatorname{Act}}
\newcommand{\rev}{\operatorname{rev}}
\newcommand{\set}[1]{\left\{#1\right\}}
\newcommand{\odd}{\operatorname{odd}}

\setlist{nosep,leftmargin=2em}
\allowdisplaybreaks

\begin{document}

\title{\bf  Characterization of Cubic Graphs with Exactly Two Lonely Edges \footnote{This research is supported by the National Natural Science Foundation of China (12301453), the Project of Scientific Research Fund of Hunan Provincial Science and Technology Department (2018WK4006).}
}
\author{Zewei Qian, Pingshan Li\footnote{ E-mail address: lips@xtu.edu.cn}
\\
{\footnotesize \em School of Mathematics and Computational Science, Xiangtan University,
Xiangtan, Hunan 411105, PR China}\\
{\footnotesize Key Laboratory of Intelligent Computing$\And$Information Processing of Ministry of Education}\\
{\footnotesize Key Laboratory for Computation and Simulation in Science and Engineering}\\
{\footnotesize National Center for Applied Mathematics in Hunan}\\
}

%\title{\bfseries Characterization of Cubic Graphs with Exactly Two Lonely Edges}
%\author{}
 \date{}
 \date{}
 \maketitle

\begin{abstract}
Perfect matchings are fundamental to understanding structural connectivity in cubic graphs, with ``lonely edges" (edges in exactly one perfect matching) forming a critical underexplored class. Goedgebeur et al. (J. Graph Theory, 112 (2026) 276-289) established an upper bound of 6 on the number of lonely edges in 3-connected non-double covered cubic graphs and completely characterized the classes $\mathcal{U}_k$ for $3\le k \le 6$. Narayana et al. subsequently gave a recursive characterization of the subclass in which the two lonely edges lie in the same perfect matching. This paper completes the classification by presenting a full structural characterization of $\mathcal{U}_2$, the family of 3-connected cubic graphs with exactly two lonely edges. We identify the minimally non-contractible ``primitive kernels" of $\mathcal{U}_2$ and develop constructive operations to generate any graph in $\mathcal{U}_2$ from these kernels. 
\end{abstract}

%\tableofcontents

\section{Introduction}\label{sec:intro}

A perfect matching is an edge set in which every vertex is incident with
exactly one edge. Perfect matchings have played a central role in graph
theory since the late nineteenth century, with applications in
combinatorial optimization, coding theory, and network design. In
particular, for cubic graphs, perfect matchings provide an important tool
for studying structural connectivity and resilience.

A fundamental theorem of Petersen states that every bridgeless cubic graph
has a perfect matching \cite{petersen1891theorie}. Schönberger subsequently
proved the corresponding edgewise statement: every edge of a bridgeless
cubic graph belongs to at least one perfect matching
\cite{schonberger1934uber}. Thus, for cubic graphs, bridgelessness is
equivalent to the property that every edge belongs to some perfect
matching. These results established the basic connection between bridges
and perfect-matching coverage.

Later research turned from the existence of perfect matchings to their
number and distribution. Esperet et al. proved that the number of perfect
matchings in a bridgeless cubic graph grows exponentially with the order of
the graph \cite{esperet2011number}. This result motivates studying the
multiplicity with which individual edges occur in perfect matchings. A
cubic graph is called \emph{matching double covered} if every edge belongs
to at least two distinct perfect matchings. Consequently, edges that occur
in exactly one perfect matching form a natural exceptional class. 
Goedgebeur et al. introduced the term \emph{lonely edge} for an edge that
belongs to exactly one perfect matching and initiated a systematic study
of non-double-covered cubic graphs \cite{GMMRW2026}. Their work builds on
an earlier result of Esperet et al., who proved that every
3-connected cubic graph is matching double covered unless it belongs to a
special class of graphs called Klee graphs
\cite{esperet2010matching}.

Klee graphs are defined recursively by starting with $K_4$ and replacing
vertices by triangles. Fiorini and Wilson conjectured that Klee graphs are precisely the
planar cubic graphs with a unique 3-edge-colouring, and a proof of this
conjecture was given in Fowler's doctoral thesis
\cite{fiorini1978edge,fowler1998thesis}. Not every Klee graph is
non-double covered. For example, the truncation of $K_4$ is still matching
double covered. Thus, a more precise structural analysis is needed to
identify exactly which Klee graphs contain lonely edges.

Goedgebeur et al. reduced the study of bridgeless non-double-covered cubic
graphs to their 3-connected components by means of 2-cut reductions
\cite{GMMRW2026}. Within the resulting 3-connected setting, they developed
a constructive framework based on vertex-to-triangle replacements and
spanning subgraph joins. For $k\geq 1$, let $\mathcal{U}_k$ denote the
subclass of graphs in this framework having exactly $k$ lonely edges. Several cases were settled in their work. They proved that no graph in the
considered class has more than six lonely edges and that
$\mathcal{U}_k=\varnothing$ for $k>6$. Moreover, they also characterized the classes with three, four, five and six lonely edges \cite{GMMRW2026}. In addition, they exhibited infinite families in \(\mathcal{U}_1\) and \(\mathcal{U}_2\), showing that the cases with at least three lonely edges and those with one or two lonely edges have substantially different structural behaviour. For example, graphs in \(\mathcal{U}_1\) and \(\mathcal{U}_2\) can contain arbitrarily many triangles, whereas those in \(\mathcal{U}_3, \mathcal{U}_4, \mathcal{U}_5, \mathcal{U}_6\) contain at most three triangles, except for \(K_4\), which contains four triangles. Recently, Narayana et al. gave a recursive characterization of the subclass in which the two lonely edges lie in the same perfect matching.

The present paper focuses on the class $\mathcal{U}_2$ of
3-connected non-double-covered cubic graphs with exactly two lonely edges.
The main objective is to give a complete structural characterization of
this class. We address the following two questions.

\begin{enumerate}
    \item Which graphs in $\mathcal{U}_2$ are minimally
    non-contractible to a smaller graph in $\mathcal{U}_2$? We call such
    graphs \emph{primitive $\mathcal{U}_2$ kernels}.

    \item Given an arbitrary graph $G\in\mathcal{U}_2$, can $G$ be
    obtained from a primitive kernel by a finite sequence of expansions,
    while every intermediate graph remains in $\mathcal{U}_2$?
\end{enumerate}

The first question is answered in Section~\ref{sec:parents}, where we
classify all primitive $\mathcal{U}_2$ kernels. The second question is
addressed in Section~\ref{sec:matching-type},
where we introduce standard markings and prove that the relevant
zero-loss expansions preserve membership in $\mathcal{U}_2$. Together,
these results provide a constructive characterization of every graph in
$\mathcal{U}_2$ from a primitive kernel.

\section{Preliminaries}\label{sec:prelim}

Throughout this paper, all graphs are assumed to be finite, simple, and cubic unless explicitly stated otherwise. Our work builds directly upon the framework established in \cite{GMMRW2026}; we therefore adopt most of the notation and conventions introduced therein, with only minor modifications made to simplify subsequent discussions. For the sake of completeness and readability, we restate key definitions below.

Let $N_G(v)=\set{x,y,z}$. When $v$ is replaced by a triangle, label the three new vertices
\(
v_x, v_y, v_z,
\)
where $v_x$ is adjacent to the old neighbor $x$, $v_y$ to $y$, and $v_z$ to $z$. Add the three triangle edges
\(
v_xv_y, v_yv_z, v_zv_x.
\)
The opposite edge of the old edge $vx$ is $v_yv_z$. This labeling depends only on the old neighbors in the current step, so every operation is completely determined at the vertex-label level.
A \emph{join rooted at \(v\)} is a spanning subgraph \(J\) of \(G\) such that \(d_J(v)=3\) and \(d_J(u)=1\) for every \(u\ne v\), where \(d_J(v)\) represents the degree of \(v\) in \(J\). Since $G$ is cubic, $J$ necessarily contains all three edges incident with $v$.
\begin{comment}
  Let $v \in V(G)$ with open neighborhood $N_G(v) = \{x_1, x_2, x_3\}$. The graph $G^v$ is obtained by deleting vertex $v$, inserting a new triangle $v_1v_2v_3$, and adding edges $v_i x_i$ for each $i = 1,2,3$. We refer to $G^v$ as the \emph{expansion of $G$ at vertex $v$}. The edge $v_jv_k$ is called the opposite edge of the old edge $vx_i$, where $\set{i,j,k}=\set{1,2,3}$.  
\end{comment}

\begin{definition}\label{def:klee-graph}
    A graph $G$ is a \textit{Klee-graph} if $G = K_4$ or there exists a Klee-graph $H$ and a vertex $w \in V(H)$ such that $G = H^w$. We denote the set of Klee-graphs by $\mathcal{K}$.
\end{definition}

Let $G\in \mathcal{K}$. Suppose that a graph $H$ can be obtained from $G$ through a finite sequence of vertex expansion operations, where each operation replaces a single vertex with a triangle. In such cases, we refer to $H$ as a \textit{descendant} of $G$. Building on this parent-descendant framework, we further define a direct hierarchical relationship: for any vertex $v\in V(G)$, the graph $G^v$ obtained by expanding $v$ into a triangle is called a \textit{child} of $G$, and conversely, $G$ is called the \textit{parent} of $G^v$. 
 Let $\PM(G)$ denote the set of perfect matchings of $G$.
An edge $e$ is called a lonely edge if it belongs to exactly one member of $\PM(G)$. Put
\[
L(G)=\set{e\in E(G):|\set{M\in\PM(G):e\in M}|=1},
\qquad \ell(G)=|L(G)|.
\]
Write
\[
\U_k=\set{G\in\K:\ell(G)=k},\qquad
\U=\bigcup_{k\ge1}\U_k.
\]

\begin{definition}\label{def:sigma}
For $G \in \mathcal{U}$ and $v \in V(G)$, define
\[
\Sigma_G(v) = \left\{ e \in L(G) : v \notin V(e), \text{ and there exists a join rooted at } v \text{ that contains } e \right\}.
\]
We simultaneously define a vertex $v$ satisfying
\(
\Sigma_G(v) = \emptyset
\)
as a \textit{zero-loss expansion vertex}, or alternatively refer to such a vertex as \textit{zero-loss expandable}.
\end{definition}

\begin{lemma}\label{lem:deletion-criterion}
Let $xy\in L(G)$ and $v\notin\set{x,y}$. Then
\(
xy\in\Sigma_G(v)
\)
if and only if
\[
\set{x,y}\cap N_G[v]=\varnothing
\quad\text{and}\quad
G-N_G[v]-\set{x,y}\text{ has a perfect matching}.
\]
\end{lemma}

\begin{proof}
If $J$ is a join rooted at $v$ containing $xy$, then $J$ contains all three edges incident with $v$. Hence every vertex in $N(v)$ already has degree $1$ in $J$, so neither $x$ nor $y$ lies in $N[v]$. After deleting $N[v]\cup\set{x,y}$, the remaining edges of $J$ are pairwise disjoint and cover all remaining vertices, and therefore form a perfect matching.

Conversely, suppose that $\set{x,y}\cap N[v]=\varnothing$ and that $M$ is a perfect matching of $G-N[v]-\set{x,y}$. Then
\(
M\cup\set{xy}\cup\set{vw:w\in N(v)}
\)
is a join rooted at $v$ containing $xy$. This proves both directions.
\end{proof}

Let $N_G(v)=\{x, y, z\}$.  Define the transfer map on lonely edges by
\[
\lambda_v(e)=
\begin{cases}
 e, & v\notin V(e),\\
 v_av_b, & e=vc {\rm~where~}\{a,b,c\}=\{x, y, z\}.
\end{cases}
\]
The following result is the set-valued form of~\cite[Lemma~3.5 and Corollary~3.6]{GMMRW2026}.

\begin{proposition}\label{prop:loss-formula}
Let $G\in\U$ and $v\in V(G)$. Then
\begin{equation}\label{eq:lonely-set-formula}
L(G^v)=\set{\lambda_v(e):e\in L(G)\setminus\Sigma_G(v)},
\end{equation}
and
\begin{equation}\label{eq:ell-formula}
\ell(G^v)=\ell(G)-|\Sigma_G(v)|.
\end{equation}
\end{proposition}

\begin{corollary}\label{cor:parent-types}
If $G\in\U_k$, then
\(
G^v\in\U_2\)  if and only if \(|\Sigma_G(v)|=k-2.
\)
\end{corollary}

\begin{corollary}\label{cor:new}
Let $G\in \mathcal{U}$ and $N_G(v)=\{x, y, z\}$. Then $\{v_xx, v_yy, v_zz\}\cap L(G^v)=\emptyset$.
\end{corollary}

\begin{proof}
Let \(w\) be an arbitrary vertex in \(\{x,y,z\}\). If \(vw\in L(G)\), then, by \cref{prop:loss-formula}, \(wv_w\notin L(G^v)\).
Otherwise, \(vw\notin L(G)\). Then \(vw\) is contained in at least two perfect matchings of \(G\). Let \(M_{vw}\) and \(M'_{vw}\) be two distinct perfect matchings of \(G\) containing \(vw\). Then \(M_{vw}-vw+wv_w+e, M'_{vw}-vw+wv_w+e\) are two distinct perfect matchings of \(G^v\), where \(e\) denotes the edge opposite \(vw\). Hence, \(wv_w\notin L(G^v)\). Combining the above, we obtain \(\{v_x x, v_y y, v_z z\}\cap L(G^v)=\emptyset.
\)
\end{proof}

\begin{definition}\label{def:primitive}
A graph $C\in\U_2$ is called a primitive $\U_2$ kernel if, for every triangle $T$ of $C$,
$C/T\notin\U_2$, where $C/T$ is the cubic graph obtained by contracting $T$ to a single vertex.
\end{definition}

\begin{definition}\label{def:neu}
For a graph family $\mathcal A\subseteq\U_2$, let $\Neu(\mathcal A)$ be the smallest family satisfying:
\begin{enumerate}[label=\textup{(\alph*)}]
\item $\mathcal A\subseteq\Neu(\mathcal A)$;
\item if $G\in\Neu(\mathcal A)$ and $\Sigma_G(v)=\varnothing$, then $G^v\in\Neu(\mathcal A)$.
\end{enumerate}
\end{definition}

\begin{proposition}\label{prop:kernel-reduction}
Let $\mathcal C$ be the set of all primitive $\U_2$ kernels. Then 
$\U_2=\Neu(\mathcal C).$

\end{proposition}

\begin{proof}
By \cref{prop:loss-formula}, expanding a graph in $\U_2$ at a vertex satisfying $\Sigma_G(v)=\varnothing$ keeps the graph in $\U_2$, so $\Neu(\mathcal C)\subseteq\U_2$. Conversely, take any $G\in\U_2$. If $G$ is not primitive, there is a triangle $T$ such that $G/T\in\U_2$. Repeating such contractions decreases the order by $2$ at every step, so the process terminates after finitely many steps at some primitive kernel $C$. Reversing the contraction sequence gives
\(
C=G_0,G_1,\cdots,G_s=G\) with \(G_{i+1}=G_i^{v_i}
\)
and  every $G_i,G_{i+1}\in\U_2$. By \cref{prop:loss-formula},  
\(
2=\ell(G_{i+1})=\ell(G_i)-|\Sigma_{G_i}(v_i)|
=2-|\Sigma_{G_i}(v_i)|.
\)
It follows that $\Sigma_{G_i}(v_i)=\varnothing$. Thus, $G\in\Neu(\mathcal C)$.
\end{proof}

Proposition~\ref{prop:kernel-reduction} shows that the only substantive task remaining is to classify the primitive kernels.

% ---- end input: prelim.tex ----
% ---- begin input: kernel_part1.tex ----
\section{The primitive kernels of $\mathcal{U}_2$}\label{sec:parents}

The triangular prism is composed of two vertex-disjoint triangles, denoted $U = u_1u_2u_3$ and $X^0 = x_1^0x_2^0x_3^0$, which are connected by the matching edges $\{{u_i x_i^0 \mid i = 1,2,3\}}$. To formalize the iterative expansion process, the concept of an ``active-end triangle'' is defined: at each step $t \geq 1$, let $X^{{t-1}} = x_1^{{t-1}}x_2^{{t-1}}x_3^{{t-1}}$ denote the triangle targeted for modification. An expansion with index $s_t \in \{{1,2,3\}}$ is performed by replacing the vertex $x_{{s_t}}^{{t-1}}$ with a new triangle $X^t = x_1^t x_2^t x_3^t$, where for each $j \neq s_t$, the new vertex $x_j^t$ is adjacent to the original vertex $x_j^{{t-1}}$. By applying a sequence of such expansions following a pattern $w = s_1 s_2 \cdots s_m$, a graph known as an \textit{{extended prism}} is constructed, which is denoted by $P(w)$.

\begin{figure}[htbp]
\centering
\includegraphics[width=0.9\textwidth]{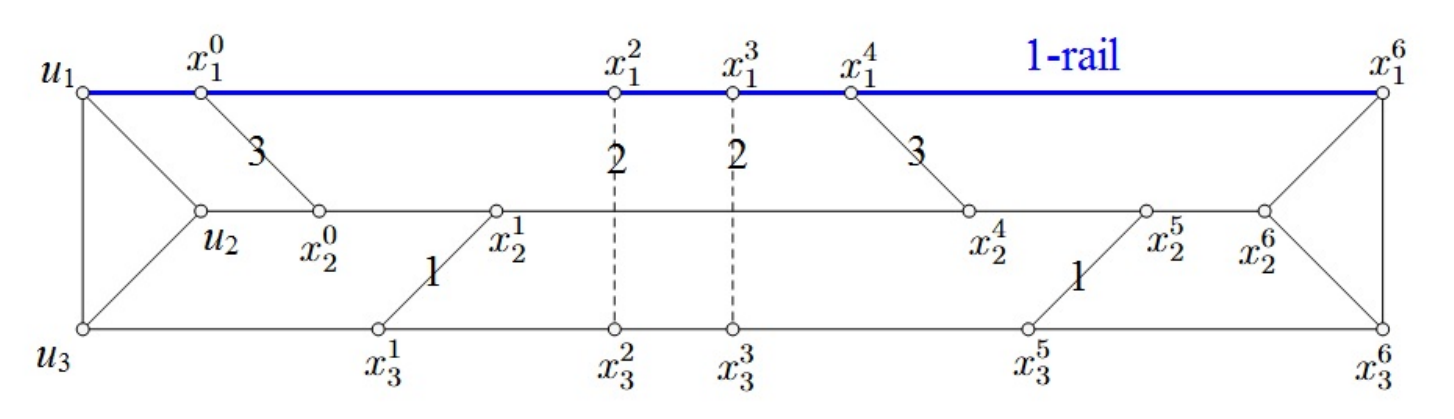}
\caption{The extended prism $P(312231)$}
\label{fig:P}
\end{figure}

It is readily observed that 
\(P(s_1s_2\cdots s_m)\) is constructed by inserting a total of \(2m\) new vertices into the edges in \(\{u_ix_i^0\mid 1\le i\le 3\}\), with the restriction that each edge is subdivided at most \(m\) times, and relabeling \(x_i^0\) as \(x_i^m\). 
The path $P_i$ with endpoints $u_i$ and $x_i^m$ obtained through this subdivision is hereafter referred to as the $i$-rail. Additionally, $m$ edges are added to connect these vertices (specifically, at step $t$, one vertex is inserted on each $i$-rail with $i\in \{1,2,3\}\setminus\{s_t\}$, and an edge is then added between these two vertices); the edge is called a \emph{$k$-rung} if both its endpoints have the superscript $k$ ($0\le k\le m-1$) and distinct subscripts. An example of $P(312231)$ is shown in Figure \ref{fig:P}.  
By flipping this graph upside down, it can be observed that $P(s_1s_2\cdots s_{m-1}s_m)$ is isomorphic to $P(s_ms_{m-1}\cdots s_2s_1)$.

Part (1) of the following lemma is implicit in the definition of extended prisms and in the proof of Theorem~3.12 in Goedgebeur et al.~\cite{GMMRW2026}, while part (2) is directly derived from the construction of extended prisms.

\begin{lemma}\label{lem:two-triangles}
{\rm (1)} A Klee graph of order at least $6$ has exactly two triangles if and only if it is an extended prism. 

{\rm (2)} For $m\ge1$, the two triangles of $P(s_1\cdots s_m)$ are precisely the fixed left end and the active right end; contracting the right end yields $P(s_1\cdots s_{m-1})$, while contracting the left end yields a graph isomorphic to $P(s_m\cdots s_2)$.
\end{lemma}

Armed with these definitions, we proceed to present the characterizations of $\mathcal{U}_i$ for $3\le i\le 6$ established by Goedgebeur et al. \cite{GMMRW2026}. 

\begin{theorem}{\rm\cite{GMMRW2026}}\label{th2} A graph $G$ belongs to $\mathcal{U}_3$ if and only if 
\begin{center}
$G\in\mathcal{W}_3\cup \{P(1^{h-1}3\,2^{p-h}): p,h\in\mathbb Z,\ p\ge 3,\ \frac{p+1}{2}\le h\le p\}$,
\end{center}
where $\mathcal{W}_3=\{T_0, T_1, T_2, T_3, P(121)\}$ and $T_0,T_1,T_2,T_3$ are shown in \cref{fig:T-branch}.
\end{theorem}

\begin{theorem}{\rm\cite{GMMRW2026}}{\label{th1}} A graph $G$ belongs to $\mathcal{U}_4$ if and only if $G\in \{P(12)\}\cup\{P(1^p): p\ge 2\}$. 
\end{theorem}

\begin{figure}[htbp]
\centering
\includegraphics[width=0.8\textwidth]{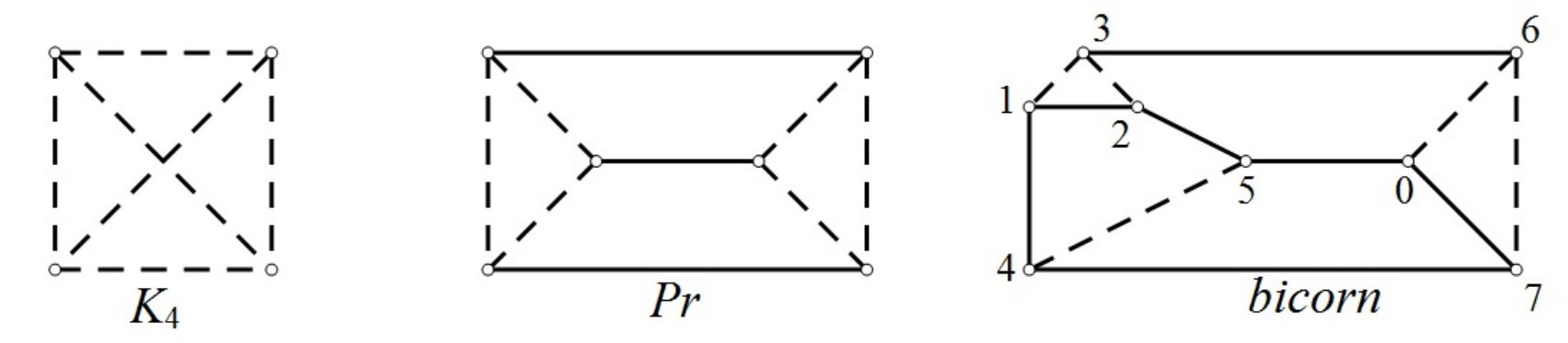}
\caption{All graphs of $\mathcal{U}_5$ and $\mathcal{U}_6$}
\label{fig:pr}
\end{figure}

\begin{theorem}{\rm\cite{GMMRW2026}}{\label{th3}} The following statements hold. 
\begin{description}
  \item[(i)] $\ell(G)\le 6$ for every $G\in \mathcal{U}$;
  \item[(ii)] The bicorn is the only element of $\mathcal{U}_5$ (see \cref{fig:pr});
  \item[(iii)] $K_4$ and the triangular prism  $Pr$ are the only elements of $\mathcal{U}_6$.
\end{description}
\end{theorem}

\subsection{The structural properties of extended prisms}

For a thorough characterization of $\mathcal{U}_2$, it is necessary to further investigate the structural properties of extended prisms. We now state several fundamental properties of such graphs. Let $w=s_1s_2\cdots s_m$ and $w'$ be words of the same length, and let $\operatorname{rev}(w)=s_ms_{m-1}\cdots s_1$. We write $w\sim w'$ if there exists a permutation $\pi\in S_3$ such that either $w'=\pi(w)$ or $w'=\pi(\operatorname{rev}(w))$, where $\pi(w)=\pi(s_1)\pi(s_2)\cdots\pi(s_m)$. Since reversal commutes with symbol permutations and applying reversal twice gives the original word, it is easy to verify that this relation is transitive.

\begin{lemma}\label{lem:word-uniqueness}
$P(w)\cong P(w')$ if and only if $w\sim w'$.
\end{lemma}

\begin{proof}
For convenience, we denote the vertices in $P(w)$ by $u_i,x_i^j$ and the vertices in $P(w')$ by $u'_i,x_i^{\prime\,j}$ (if they exist), where $i\in\{1,2,3\}$ and $j\in\{0,1,\ldots,m\}$. Let $w=s_1s_2\cdots s_m$ and $w'=s'_1s'_2\cdots s'_m$. We prove the two directions separately.

\medskip
\noindent\textbf{Sufficiency.}
Suppose that $w\sim w'$. Then there exists $\pi\in S_3$ such that $w'=\pi(w)$ or $w'=\pi(\operatorname{rev}(w))$. For $w'=\pi(w)$, we construct an isomorphism between $P(w)$ and $P(w')$ as follows:
\[
u_i\longmapsto u'_{\pi(i)},\qquad
x_i^t\longmapsto x_{\pi(i)}^{\prime\,t},
\]
where $0\le t\le m$, and the map is defined only for vertices that exist. By the construction of extended prisms, this bijection preserves adjacency.

For $w'=\pi(\operatorname{rev}(w))$, we first show that $P(w)\cong P(\operatorname{rev}(w))$. Denote the vertices of $P(\operatorname{rev}(w))$ by $\widetilde u_i,\widetilde x_i^t$. Define
\[
u_i\longmapsto \widetilde x_i^m,\qquad
x_i^m\longmapsto \widetilde u_i,\qquad
x_i^t\longmapsto \widetilde x_i^{m-t-1},
\]
where $0\le t\le m-1$ in the last formula, again considering only existing vertices. Indeed, $x_i^t$ exists precisely when $i\ne s_{t+1}$, which is also the condition for $\widetilde x_i^{m-t-1}$ to exist. This bijection exchanges the two end triangles, reverses each rail, and maps each rung to the corresponding rung. Hence it is an isomorphism. By the first case, $P(\operatorname{rev}(w))\cong P(\pi(\operatorname{rev}(w)))$. Thus $P(w)\cong P(w')$.

\medskip
\noindent\textbf{Necessity.}
Suppose that $\varphi$ is an isomorphism between $P(w)$ and $P(w')$. Since every extended prism has exactly two triangles, namely its two end triangles, $\varphi$ must either preserve these two triangles individually, in which case
\[
\varphi(u_1)=u'_{i_1},\qquad
\varphi(u_2)=u'_{i_2},\qquad
\varphi(u_3)=u'_{i_3},
\]
or interchange them, in which case
\[
\varphi(u_1)=x_{i_1}^{\prime\,m},\qquad
\varphi(u_2)=x_{i_2}^{\prime\,m},\qquad
\varphi(u_3)=x_{i_3}^{\prime\,m}.
\]

For the first case, let
\(
\pi=
\begin{pmatrix}
1&2&3\\
i_1&i_2&i_3
\end{pmatrix}.
\)
We prove by induction on $m$ that $w'=\pi(w)$ and that $\varphi(x_i^m)=x_{\pi(i)}^{\prime\,m}$ for every $i\in\{1,2,3\}$. When $m=0$, both graphs are triangular prisms, and the three edges joining their end triangles force the same permutation on both ends. Suppose that $m\ge1$. Contracting the right end triangle in each graph induces an isomorphism
\[
\overline{\varphi}:
P(s_1\cdots s_{m-1})
\longrightarrow
P(s'_1\cdots s'_{m-1}).
\]
This induced isomorphism preserves the left end and hence also the right end, since each contracted graph has exactly two triangles. By the induction hypothesis,
\[
s'_t=\pi(s_t)\quad(1\le t<m),
\qquad
\overline{\varphi}(x_i^{m-1})
=x_{\pi(i)}^{\prime\,m-1}\quad(i\in\{1,2,3\}).
\]
The contracted vertices are $x_{s_m}^{m-1}$ and $x_{s'_m}^{\prime\,m-1}$, respectively. Therefore,
\(
x_{\pi(s_m)}^{\prime\,m-1}
=\overline{\varphi}(x_{s_m}^{m-1})
=x_{s'_m}^{\prime\,m-1},
\)
which gives $s'_m=\pi(s_m)$ and hence $w'=\pi(w)$. To complete the inductive assertion, observe that, for $j\ne s_m$, the vertex $x_j^{m-1}$ survives the expansion and $x_j^m$ is its unique neighbor in the new triangle. Thus $\varphi(x_j^m)=x_{\pi(j)}^{\prime\,m}$. The remaining vertex of the new triangle must also map to the vertex with the corresponding index. This completes the induction and proves $w\sim w'$ in the first case.

For the second case, let
\[
\rho:P(w')\longrightarrow P(\operatorname{rev}(w'))
\]
be the end-interchanging isomorphism constructed in the sufficiency proof. Then $\rho\circ\varphi$ preserves the two ends individually. By the first case, there exists $\pi\in S_3$ such that $\operatorname{rev}(w')=\pi(w)$. Since reversal commutes with symbol permutations, we obtain $w'=\pi(\operatorname{rev}(w))$. Hence $w\sim w'$, completing the proof.
\end{proof}

We decompose each word $w$ into a product of maximal constant blocks:
\[
w = a_1^{r_1}a_2^{r_2}\cdots a_q^{r_q},
\qquad a_i\ne a_{i+1}\ (1\le i<q),\qquad r_i\ge1\ (1\le i\le q).
\]
The integer $q$ is called the \emph{number of runs} of $w$. Combining \cref{th2,th1} with \cref{lem:word-uniqueness} gives a formulation independent of the choice of normalization.
\begin{lemma}\label{lem:H-language}
Let $|w|\ge4$. Then $P(w)\in\U_3\cup\U_4$ if and only if the run decomposition of $w$ satisfies one of the following three conditions:
\begin{enumerate}[label=\textup{(\alph*)}]
\item $q=1$;
\item $q=2$, and at least one of the two end runs has length $1$;
\item $q=3$, the middle run has length $1$, and $a_1\ne a_3$.
\end{enumerate}
\end{lemma}

\begin{proof}
By \cref{th1} and \cref{lem:word-uniqueness}, we have that constant words give the $\U_4$ family. By \cref{lem:word-uniqueness}, a two-run word with a one-letter end run can, after reversal and a global relabeling, be written as $1^a3$ with $a\ge3$, since $|w|\ge4$. By \cref{th2}, $P(w)\in\U_3$. In the three-run case, after reversal we may assume that the first block is at least as long as the third, and then globally relabel it as $1^a3\,2^b$. Recall that $|w|\ge 4$. By \cref{th2}, $P(w)\in \U_3$.

Conversely, suppose that $P(w)\in\U_3\cup\U_4$.
Since $|w|\ge4$, the graph $P(w)$ has order at least $14$.
If $P(w)\in\U_4$, then \cref{th1} excludes the exceptional graph
$P(12)$ and gives $P(w)\cong P(1^p)$ for some $p\ge2$.
By \cref{lem:word-uniqueness}, $w\sim1^p$, so $w$ is constant
and satisfies~(a).

Now suppose that $P(w)\in\U_3$.
By \cref{lem:two-triangles}, $P(w)$ has exactly two triangles,
whereas each of $T_0,T_1,T_2,T_3$ has three triangles.
Moreover, $P(w)\not\cong P(121)$, since $P(121)$ has order $12$.
Thus \cref{th2,lem:word-uniqueness} give
\(
w\sim1^{h-1}3\,2^{p-h}
\)
for some integers $p,h$ satisfying $p=|w|\ge4$ and
$(p+1)/2\le h\le p$.
If $h=p$, this word has two runs, with the last run of length $1$,
so it satisfies~(b). If $h<p$, it has three runs, the middle run
has length $1$, and the first and last symbols are distinct,
so it satisfies~(c).
Since conditions~(a)--(c) are preserved under reversal and global
relabeling, the corresponding condition also holds for $w$.
This proves necessity.
\end{proof}

%; the middle block must consist of a single letter, and the two outer symbols must be different. The converse is identical. The assumption $|w|\ge4$ excludes the additional exceptional short patterns $12$ and $121$.

As a conclusion to this subsection, a result concerning the parent graphs of $\mathcal{U}_2$ primitive kernels is presented, which plays a crucial role in the classification of all $\mathcal{U}_2$ primitive kernels in the subsequent two subsections.

\begin{lemma}\label{lem:primitive-parent}

Let $C$ be a primitive $\mathcal{U}_2$ kernel, let $T$ be any triangle of $C$, and put $H = C/T$. Then
\(\
H \in \mathcal{U}_3 \cup \mathcal{U}_4.
\)
\end{lemma}

\begin{proof}
Write $C = H^v$. Because $C\in \mathcal{U}_2$, $C$ is a Klee graph. By \cite[Lemma~3.2]{GMMRW2026}, $H$ is also a Klee graph.
If $H\notin \mathcal{U}$, then every edge of $H^v$ in at least two perfect matchings:
for an edge of the new triangle, use the matchings containing the
opposite old edge, and for any other edge use the matchings containing
its corresponding old edge. This contradicts $C\in\U_2$.
Hence $H\in\U$. By \cref{prop:loss-formula},
\(
2 = \ell(C) = \ell(H) - |\Sigma_H(v)|.
\)
It follows that $\ell(H) \geq 2$. 
Since \(C\) is a primitive \(\U_2\) kernel and \(H = C/T\), we have \(\ell(H) \ge 3\). Note that Klee graphs are obtained from \(K_4\) by expansions, and expansion never increases the number of lonely edges, so \(\ell(H) \leq 6\). To complete the proof, it suffices to show that neither \(\ell(H)=5\) nor \(\ell(H)=6\).

By \cref{th3}, $\mathcal{U}_6$ consists only of $K_4$ and the triangular prism $\text{Pr}$, while $\mathcal{U}_5$ consists only of the bicorn (see \cref{fig:pr}). 
Thus, we list all their children up to isomorphism in ~\cref{tab:descendants}, which shows that none of their children contains exactly two lonely edges.
\begin{table}[ht]
\centering
\caption{Children of graphs in $\mathcal{U}_6$ and $\mathcal{U}_5$.}
\begin{tabular}{cccc}
\toprule
Parent graph & Expanded vertex & Child graph & Lonely edges \\
\midrule
$K_4$ & all vertices & $\text{Pr}$ & 6 \\
$\text{Pr}$ & all vertices & bicorn & 5 \\
\multirow{3}{*}{bicorn} & $\{0,1,2,7\}$ & $P(12)$ & 4 \\
& $\{3,6\}$ & $P(11)$ & 4 \\
& $\{4,5\}$ & $T_0$ & 3 \\
\bottomrule
\end{tabular}
\label{tab:descendants}
\end{table}

Thus, $\ell(H)$ cannot be $5$ or $6$.
\end{proof}

\begin{comment}
    To further restrict possible values of $\ell(H)$, note that the above classification exhausts all graphs in $\mathcal{U}_5$ and $\mathcal{U}_6$, and Table \ref{tab:descendants} lists all their nonisomorphic children. The number of lonely edges of a descendant can only be 6, 5, 4, or 3, and is never 2. This implies there is no one-step expansion $\mathcal{U}_5 \to \mathcal{U}_2$ or $\mathcal{U}_6 \to \mathcal{U}_2$, so $\ell(H)$ cannot be 5 or 6. 
    \end{comment}

Primitive $\mathcal{U}_2$ kernels are divided into two distinct classes based on whether they contain exactly two triangles: those with exactly two triangles are referred to as \textit{extended-prism kernels}, while the others are termed \textit{non-extended-prism kernels}. The subsequent two sections are dedicated to the characterization of these two classes of primitive $\mathcal{U}_2$ kernels, respectively.

\subsection{Primitive extended-prism  $\mathcal{U}_2$ kernels}

The primary objective of this subsection is to present a characterization of all primitive kernels in $\mathcal{U}_2$ that belong to the class of extended prisms.

%This proof works directly at the level of runs. Hence compatibility with global symbol permutations and reversal is already built into \cref{lem:word-uniqueness,lem:H-language}, and no separate, potentially inconsistent normalizations of the two parent patterns are needed.
%\subsection{The complete finite table at length four}

\begin{lemma}\label{lem:length4-kernels}
The primitive extended-prism $\U_2$ kernels of length four are exactly
\[
P(1121),P(1122),P(1212),P(1213),P(1221),P(1223),P(1231).
\]
\end{lemma}

\begin{proof}

Up to global relabeling and reversal, there are exactly ten patterns of length four:
\[
1111,\quad1112,\quad 1121,\quad1122,\quad1123,\quad1212,\quad 1213,\quad1221,\quad1223,\quad1231.
\]
To establish the validity of the conclusion, it suffices to identify all $\mathcal{U}_2$ graphs among these 10 extended prisms, and then prove that exactly 7 of these graphs satisfy the property that contracting any triangle results in a graph no longer belonging to $\mathcal{U}_2$. By exactly enumerating all perfect matchings of the aforementioned graphs (this can be verified manually: for each selected edge, find the perfect matchings that contain it; once an edge is contained in two distinct perfect matchings, we no longer need to find all perfect matchings containing this edge. See also the Python code in~\cite{qianli}), we can obtain the results presented in \Cref{tab:length4}.

\begin{table}[htbp]
\centering
\caption{Length-four expansion patterns ($m(G):$ number of perfect matchings).}
\label{tab:length4}
\begin{tabular}{c r r c@{\qquad}c r r c}
\toprule
Pattern & $m(G)$ & $\ell(G)$ & Class & Pattern & $m(G)$ & $\ell(G)$ & Class\\
\midrule
1111&15&4&$\U_4$ & 1112&11&3&$\U_3$\\
1121&9&2&$\U_2$ & 1122&11&2&$\U_2$\\
1123&11&3&$\U_3$ & 1212&9&2&$\U_2$\\
1213&9&2&$\U_2$ & 1221&10&2&$\U_2$\\
1223&9&2&$\U_2$ & 1231&10&2&$\U_2$\\
\bottomrule
\end{tabular}
\end{table}

Note that both end contractions of each such graph yield length-three patterns, and $111,112,121, 123$ represent all length-three patterns up to global relabeling and reversal. Since all four length-three equivalence classes belong to $\U_3\cup\U_4$, all seven graphs in $\mathcal{U}_2$ are primitive.
\end{proof}

% ---- end input: kernel_part1.tex ----
%\subsection{Ladder-chain construction for the infinite candidates}\label{subsec:ladder-kernels}

\begin{comment} 
Let
\(
A_r=P(12^r), r\ge2.
\)
Use the natural labeling from Section~\ref{sec:parents}: the fixed left-end triangle is $U=u_1u_2u_3$, the initial right end is $X^0=x_1^0x_2^0x_3^0$, and the final active end is
$X^{r+1}=x_1^{r+1}x_2^{r+1}x_3^{r+1}$. \end{comment}
For $e\in L(G)$, let $M_e$ denote the unique perfect matching containing $e$. Put
\[
\beta_r=
\begin{cases}
 x_1^{r+1}x_2^{r+1},&r\text{ even},\\
 x_2^{r+1}x_3^{r+1},&r\text{ odd}.
\end{cases}
\]

\begin{lemma}\label{lem:ladder-lonely}
For every $r\ge2$,
\(
L(P(12^r))=\set{u_1u_3, x_2^0x_3^0,\beta_r}.
\)
\begin{comment}
 Moreover   
 the unique perfect matching containing $u_1u_3$ is the same as the unique perfect matching containing $\beta_r$; denote it by $M_r$. Let $N_r$ be the unique perfect matching containing $x_2^0x_3^0$. Then $M_r\ne N_r$. In addition, \end{comment}
 Moreover, $M_{u_1u_3}=M_{\beta_r}\not=M_{x_2^0x_3^0}$ and the three vertices of the active end satisfy
\begin{align}
\Sigma_{P(12^r)}(x_1^{r+1})&=
\begin{cases}\set{x_2^0x_3^0},&r\text{ even},\\\set{u_1u_3},&r\text{ odd},\end{cases}
\label{eq:ladder-left-loss}\\
\Sigma_{P(12^r)}(x_2^{r+1})&=\varnothing,\label{eq:ladder-middle-loss}\\
\Sigma_{P(12^r)}(x_3^{r+1})&=
\begin{cases}\set{u_1u_3},&r\text{ even},\\\set{x_2^0x_3^0},&r\text{ odd}.\end{cases}
\label{eq:ladder-right-loss}
\end{align}
\end{lemma}

\begin{proof}
By \cref{lem:word-uniqueness}, we have $P(12^r)\cong P(2^r1)\cong P(1^r32^0)$. 
By \cref{th2}, $P(12^r)\in\U_3$, so it is enough to identify three lonely edges.

First fix $u_1u_3$. The vertex $u_2$ is forced to match $x_2^0$, which in turn forces $x_3^0$ to match $x_3^1$. Continuing along the two ladder tracks forces the matching edges successively to be
\[
x_1^1x_1^2,
\ x_3^2x_3^3,
\ x_1^3x_1^4,
\ x_3^4x_3^5,\cdots,
\]
and the last two unmatched vertices at the active end are joined precisely by $\beta_r$. Thus the perfect matching containing $u_1u_3$ is unique and also contains $\beta_r$. Running the same forcing argument backward from the right end shows that the perfect matching containing $\beta_r$ must be this same matching. Thus, $M_{u_1u_3}=M_{\beta_r}$.

If instead $x_2^0x_3^0$ is fixed, then the left end is forced to use $u_2u_3$ and $u_1x_1^1$. The forced track then becomes
\[
x_3^1x_3^2,
\ x_1^2x_1^3,
\ x_3^3x_3^4,
\ x_1^4x_1^5,\cdots,
\]
and terminates at the active end with the edge between the two remaining vertices. Hence $M_{x_2^0x_3^0}$ is also unique; it does not contain $u_1u_3$ at the left end, so it is different from $M_{u_1u_3}$. Since $P(12^r)\in\U_3$, these three edges are all the lonely edges.

Finally consider joins rooted at a vertex of the active end.
First let the root be $v=x_j^{r+1}$, where $j\in\{1,3\}$, and put
$q=4-j$. The three edges incident with $v$ are selected. In particular,
the other two vertices of the active-end triangle already have degree
$1$, so $x_2^0x_2^{r+1}$ and $x_q^rx_q^{r+1}$ are excluded.
The edge $\beta_r$ either is incident with $v$ or joins the other two
vertices of the triangle; in the latter case it cannot be selected.
Thus $\beta_r\notin\Sigma_{P(12^r)}(v)$.

The selected rail edge $x_j^rx_j^{r+1}$ determines the rail choices
backward through the internal rungs. Indeed, at each rung
$x_1^kx_3^k$ ($1\le k\le r$), exactly one selected edge enters from
the right. Its endpoint already has degree $1$, so the rung is
excluded and the other endpoint must use its left rail edge.
Consequently, the selected rail switches at each rung. At the left
boundary, the selected edge is $u_1x_1^1$ precisely when
\(j=1\) and \( r\) is even or
\(j=3\) and \(r\) is odd;  
in the other two cases it is $x_3^0x_3^1$.

If $u_1x_1^1$ is selected, then $u_1u_3$ is excluded.
Since $x_2^0x_2^{r+1}$ and $x_3^0x_3^1$ are also excluded,
the remaining four vertices at the left end can be covered by either
\(
\{u_2u_3,x_2^0x_3^0\}
\) or \(
\{u_2x_2^0,u_3x_3^0\}.
\)
The first choice gives a join containing $x_2^0x_3^0$.
If $x_3^0x_3^1$ is selected instead, then $x_2^0x_3^0$ is excluded,
and the degree requirements force $u_2x_2^0$ and $u_1u_3$.
In each case, the indicated left-end edges, the alternating rail
edges and the three edges incident with $v$ give a spanning subgraph
in which $v$ has degree $3$ and every other vertex has degree $1$.
This proves both the existence of a join containing the stated lonely
edge and the exclusion of the other nonincident lonely edges.
Hence~\eqref{eq:ladder-left-loss} and~\eqref{eq:ladder-right-loss} hold.

For the middle root $x_2^{r+1}$, every join contains
$x_2^0x_2^{r+1}$, so it cannot contain $x_2^0x_3^0$.
It cannot contain $u_1u_3$ either, since that would force $u_2x_2^0$
and give the non-root $x_2^0$ degree $2$.
The remaining lonely edge $\beta_r$ is incident with the root and
therefore does not belong to its loss set.
This proves~\eqref{eq:ladder-middle-loss}.
\end{proof}

Define
\(
\U_2^{=} = \set{G\in\U_2:M_\alpha=M_\beta},
\U_2^{\ne} = \set{G\in\U_2:M_\alpha\ne M_\beta}
\)
where $L(G)=\set{\alpha,\beta}$. Let $M$ be a perfect matching of $G$ containing $vx$, where $N_G(v)=\{x,y,z\}$. Define its natural lift by
\begin{equation}\label{eq:matching-lift}
\widehat M^{\,v}
=(M\setminus\set{vx})\cup\set{v_xx,v_yv_z}.
\end{equation}

\begin{corollary}\label{cor:infinite-kernel-type}
For every $r\ge2$, both $P(12^r1)$ and $P(12^r3)$ belong to $\U_2$. Their types are
\[
\begin{array}{c|cc}
& P(12^r1)&P(12^r3)\\ \hline
r\text{ even}&\U_2^{=}&\U_2^{\ne}\\
r\text{ odd}&\U_2^{\ne}&\U_2^{=}
\end{array}
\]
\end{corollary}

\begin{proof}
The graph $P(12^r1)$ is obtained by expanding $P(12^r)$ at $x_1^{r+1}$, and $P(12^r3)$ by expanding at $x_3^{r+1}$. By \cref{prop:loss-formula,lem:ladder-lonely}, each expansion loses exactly one lonely edge, so the resulting graph lies in $\U_2$. The unique perfect matching of each surviving lonely edge is the natural lift of the corresponding original unique perfect matching. Combining $M_{u_1u_3}=M_{\beta_r}\ne M_{x_2^0x_3^0}$ with~\eqref{eq:ladder-left-loss} and~\eqref{eq:ladder-right-loss} gives the parity classification in the table immediately.
\end{proof}

\begin{theorem}\label{thm:extended-kernels}
The primitive extended-prism $\U_2$ kernels are exactly
\[
\mathcal E_2=
\set{P(1121),P(1122),P(1212),P(1213),P(1231)}
\cup\set{P(12^r1),P(12^r3):r\ge2}.
\]
\end{theorem}

\begin{proof}
First we show that all listed graphs are primitive. The five isolated patterns and the two cases $P(1221),P(1223)$ with $r=2$ follow from \cref{lem:length4-kernels}. Let $r\ge3$. By \cref{cor:infinite-kernel-type}, $P(12^r1),P(12^r3)\in\U_2$. Each has only the two end triangles. Contracting the right end gives $P(12^r)\in\U_3$; after contracting the left end, the resulting pattern is $t2^r$ with $t=1$ or $3$, which after a global relabeling is also equivalent to $12^r$. Thus this parent also belongs to $\U_3$. Hence neither graph has a $\U_2$ parent, and both are primitive.

For completeness, let $C=P(w)$ be a primitive extended-prism $\U_2$ kernel, with $m=|w|$. By the proof of \cref{lem:length4-kernels}, extended prisms of length at most three lie in $\U_6,\U_5,\U_4$, or $\U_3$. It follows that $m\ge4$. If $m=4$, \cref{lem:length4-kernels} gives seven patterns, exactly the five isolated patterns listed in the theorem together with the two family members at $r=2$.

Now suppose $m\ge5$. By primitivity and \cref{lem:primitive-parent}, both parent graphs obtained by contracting an end triangle belong to $\U_3\cup\U_4$. By \cref{lem:two-triangles}(2), this is equivalent to stating that both the length-$(m-1)$ prefix and suffix of the corresponding word satisfy the conditions specified in \cref{lem:H-language}. Notably, the full word $w$ cannot satisfy this language; if it did, $C$ would belong to $\U_3\cup\U_4$, which contradicts the assumption that $C\in\U_2$. With this in mind, we claim that 
\(
w\sim12^{m-2}1\) or \(w\sim12^{m-2}3.
\)

To justify this claim, we begin by considering the run decomposition of $w$:
\(
w=a_1^{r_1}\cdots a_q^{r_q}.
\)
Deleting either the first or last letter of $w$ removes at most one terminal run from this decomposition. Given that both the prefix $w^-=a_1^{r_1}\cdots a_{q-1}^{r_{q-1}}a_q^{r_q-1}$ and suffix $w^+=a_1^{r_1-1}a_2^{r_2}\cdots a_q^{r_q}$ contain at most three runs, it immediately follows that $q\le4$. Furthermore, since $w$ itself does not satisfy the language defined in \cref{lem:H-language}, we can conclude that $q\ge2$.

\begin{itemize}
    \item Suppose $q=4$. In order for deletion of either the first or last letter to leave at most three runs, we must have $r_1=r_4=1$. The three run lengths of $w^-$ are then $(1,r_2,r_3)$, and \cref{lem:H-language}(c) gives $r_2=1$; similarly $w^+$ gives $r_3=1$. Thus $m=4$, contradicting $m\ge5$.
    \item Suppose $q=2$. Note that $w$ does not satisfy \cref{lem:H-language}(b). We find that $r_1,r_2\ge2$. The run lengths of the suffix are $(r_1-1,r_2)$, and one of its end runs must have length $1$, so $r_1=2$. Similarly the prefix forces $r_2=2$, again giving $m=4$, a contradiction.
    \item Suppose $q=3$. If $r_1>1$ and $r_3>1$, then both the prefix and suffix still have three runs, and both force $r_2=1$ and $a_1\ne a_3$, so $w$ itself satisfies \cref{lem:H-language}(c), a contradiction. Hence at least one end run has length $1$. Without loss of generality, let $r_1=1$. If $r_3>1$, then the prefix still has three runs and again forces $r_2=1$ and $a_1\ne a_3$, so $w$ satisfies the condition, a contradiction. Thus $r_3=1$. Symmetrically, $r_3=1$ also forces $r_1=1$. Therefore,
\(w=a_1a_2^{m-2}a_3, a_1\ne a_2, a_3\ne a_2.\)
 If $a_1=a_3$, a global relabeling gives $12^{m-2}1$; if $a_1\ne a_3$, the three symbols are distinct, and a global relabeling gives $12^{m-2}3$.
\end{itemize}
Putting $r=m-2\ge3$ gives the two infinite families. This proves completeness.
\end{proof}
% ---- begin input: kernel_part2.tex ----
\subsection{Primitive non-extended-prism kernels}\label{sec:finite-branch}

In this subsection, a characterization of primitive non-extended-prism $\mathcal{U}_2$ kernels is presented. 

\begin{comment}
   To make the vertex-by-vertex certificates reproducible, fix the following labeled expansion. Suppose $H$ has vertex set $\set{0,1,\cdots,n-1}$ and
\(
N_H(v)=\set{a,b,c \mid a<b<c}.
\)
Define $\mathsf B_v(H)$ as follows: delete the edges $va,vb$ and keep $vc$; add new vertices $n,n+1$ and the edges
\(
v n, v(n+1),n(n+1), na,(n+1)b.
\)

The new triangle is $\set{v,n,n+1}$. Starting from the labeled $K_4$ with vertices $0,1,2,3$, define
\begin{equation}\label{eq:T-def}
T_0=\mathsf B_4\bigl(\mathsf B_0(\mathsf B_0(K_4))\bigr),
T_1=\mathsf B_0(T_0),
T_2=\mathsf B_2(T_1),
T_3=\mathsf B_9(T_2),
\end{equation}
and
\begin{equation}\label{eq:D-def}
D_0=\mathsf B_1(T_0),\qquad
D_1=\mathsf B_{10}(T_1),\qquad
D_2=\mathsf B_0(T_1).
\end{equation} 
\end{comment}

\begin{figure}[htbp]
\centering
\includegraphics[width=1\textwidth]{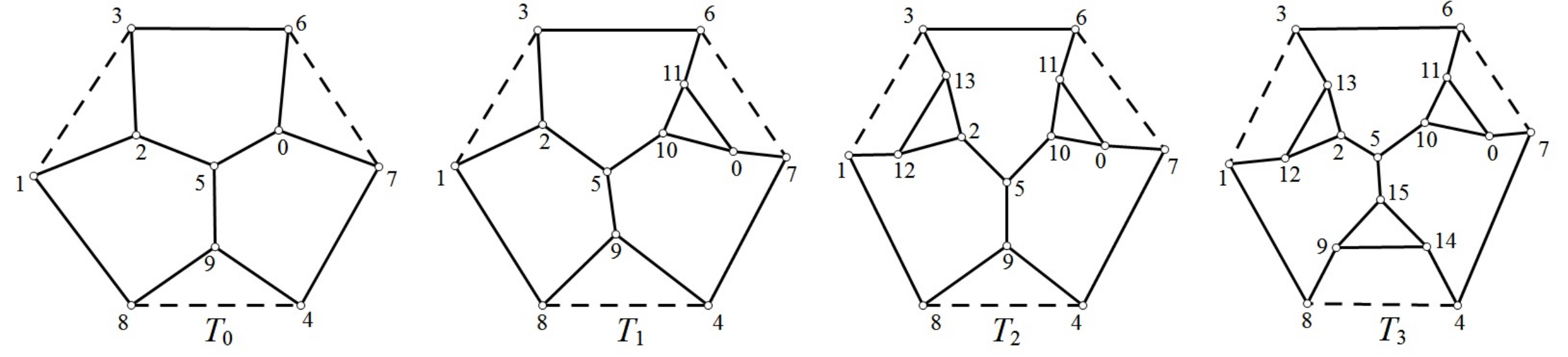}
\caption{The finite tricorn branch. In all four graphs the lonely edges are $(1,3),(4,8),(6,7)$.}
\label{fig:T-branch}
\end{figure}

\begin{figure}[htbp]
\centering
\includegraphics[width=0.78\textwidth]{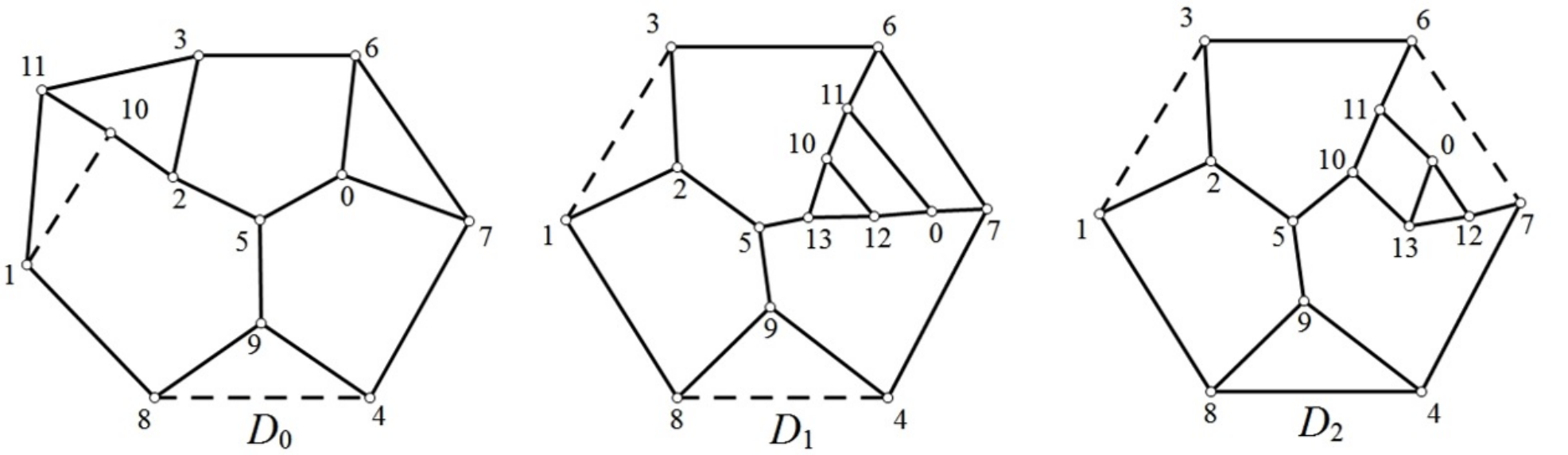}
\caption{The three primitive kernels that are not extended prisms.}
\label{fig:D-kernels}
\end{figure}

Let $D_0, D_1, D_2, T_0, T_1, T_2, T_3$ be the graphs shown in \cref{fig:T-branch,fig:D-kernels}. Dashed edges in the figures are lonely edges. Direct computation gives
\begin{equation}\label{eq:T-lonely}
L(T_i)=\set{(1,3),(4,8),(6,7)}\qquad(0\le i\le3),
\end{equation}
and
\begin{equation}\label{eq:D-lonely}
L(D_0)=\set{(1,10),(4,8)},L(D_1)=\set{(1,3),(4,8)}, L(D_2)=\set{(1,3),(6,7)}.
\end{equation}

\begin{lemma}\label{lem:finite-certificate}
The set consisting of all children of $T_0, T_1, T_2, T_3$ that are primitive $\mathcal{U}_2$ kernels is $\{D_0, D_1, D_2\}$.
\end{lemma}

\begin{proof}

The idea of the proof is relatively simple, as follows. For each vertex $v_{ij}\in V(T_i)$, where $i\in\{0,1,2,3\}$ and $1\leq j\leq |V(T_i)|$, we first enumerate the perfect matchings of $T_i^{v_{ij}}$ to determine its number of lonely edges. For this purpose, the number of perfect matchings containing each edge only needs to be counted up to two. We retain only those graphs having exactly two lonely edges. For each retained graph, we then contract each triangle separately and determine the number of lonely edges in the resulting graph. If any such contraction yields a graph with exactly two lonely edges, the original graph is discarded; otherwise, it is a primitive $\mathcal{U}_2$ kernel. Finally, we remove isomorphic duplicates, retaining one representative from each isomorphism class. Since this verification involves numerous cases, we provide an algorithmic verification in~\cite{qianli}.
\end{proof}

\begin{lemma}{\rm\cite{GMMRW2026}}\label{lem4} A graph $G\in\mathcal U_3$ has exactly three triangles
if and only if $G$ is isomorphic to one of
$T_0,T_1,T_2,T_3$. Moreover, every Klee graph of order at least $10$ with exactly three triangles is isomorphic to $T_0$ or to a descendant of $T_0$. 
\end{lemma}

\begin{theorem}\label{thm:nonextended-kernels}
The primitive $\U_2$ kernels that are not extended prisms are exactly
\(
D_0,D_1,D_2.
\)
\end{theorem}

\begin{proof}
Note that each graph in $\{D_0, D_1,D_2\}$ has three triangles and therefore is not an extended prism. By \cref{lem:finite-certificate}, \(
D_0,D_1,D_2
\) are primitive non-extended-prism $\U_2$ kernels. 
\begin{comment}
 By \cref{lem:finite-certificate}(b)(c), $D_0,D_1,D_2\in\U_2$, and every triangle-contraction parent of each belongs to $\U_3\cup\U_4$. Hence none has a $\U_2$ parent, so all three are primitive. Each has three triangles and therefore is not an extended prism.   
\end{comment}

Conversely, let $C$ be a primitive $\U_2$ kernel that is not an extended prism. By \cref{lem:two-triangles}, $C$ has at least three triangles. By \cref{lem:primitive-parent}, every triangle-contraction parent of $C$ belongs to $\U_3\cup\U_4$. According to~\cref{th2,th1}, graphs in $\U_4$ have exactly two triangles and graphs in $\U_3$ have at most three. Therefore $C$ has at most four triangles: contracting one triangle reduces the number of triangles by at most one, while the parent has at most three.

If $C$ has exactly three triangles, then by \cref{lem4} it is a descendant of $T_0$. Choose a shortest expansion sequence from $T_0$ to $C$, and let $H$ be the penultimate graph. Then $H$ is a triangle-contraction parent of $C$, so $H\in\U_3\cup\U_4$. As a descendant of $T_0$, the graph $H$ has at least three triangles, so $H\notin\U_4$ and hence $H\in\U_3$. By \cref{lem4}, $H\in\set{T_0,T_1,T_2,T_3}$, and therefore $C$ is a child of some $T_i$. By \cref{lem:finite-certificate}, $C\in \{D_0, D_1, D_2\}$.

 If $C$ has exactly four triangles, then contracting any triangle gives a parent with at least three triangles. Since the parent has at most three, it is a $\U_3$ graph with exactly three triangles and hence again belongs to $\set{T_0,T_1,T_2,T_3}$. Thus $C$ is again a child of some $T_i$. By \cref{lem:finite-certificate}, $C\in \{D_0, D_1, D_2\}$. This contradicts the fact that each $D_i$ has exactly three triangles. Hence this case cannot occur.
%Thus, every candidate appears in the complete vertex-by-vertex table of \cref{lem:finite-certificate}(b). That table shows that the primitive ones are exactly $D_0,D_1,D_2$; every other descendant has another $\U_2$ parent. This proves completeness.
\end{proof}

% ---- end input: kernel_part2.tex ----

\section{Characterization of $\mathcal{U}_2$}\label{sec:matching-type}

%We now address how expansion proceeds after the kernel stage. The structural criterion uses only the current graph, its two lonely edges, and their unique perfect matchings; the full expansion history from $K_4$ is no longer recorded.

In this section, we will provide a complete characterization of $\mathcal{U}_2$ using mathematical induction. As preliminary preparation, we first introduce two lemmas. %which include two types of markings and three operation modes.

\begin{lemma}\label{lem:matching-type-preserved}
Suppose \(L(G)=\set{\alpha,\beta}\) and \(G^v\in\U_2\). %Then the unique perfect matchings of \(\lambda_v(\alpha)\) and \(\lambda_v(\beta)\) are \(\widehat M_\alpha^{\,v}\) and \(\widehat M_\beta^{\,v}\), respectively, as defined in \cref{eq:matching-lift}. In particular,
Then \(
M_\alpha=M_\beta
\Longleftrightarrow
M_{\lambda_v(\alpha)}=M_{\lambda_v(\beta)}.
\)
\end{lemma}

\begin{proof}
By \cref{prop:loss-formula}, after a zero-loss expansion the two lonely edges are exactly $\lambda_v(\alpha)$ and $\lambda_v(\beta)$. Contracting the new triangle in a perfect matching containing either transferred edge gives a perfect matching containing the corresponding original lonely edge. The original lonely edge belongs to only one perfect matching, and~\eqref{eq:matching-lift} gives its unique lift. The equivalence follows immediately.
\end{proof}

\begin{lemma}\label{lem:new-triangle-active}
Let \( v \) be a zero-loss expansion vertex of \( G \) where \( G \in \mathcal{U}_2 \). For \( x \in N_G(v) \), the new vertex \( v_x \) is zero-loss expandable in \( G^v \) if and only if, for every lonely edge \( e \) not incident with \( v \), \(vx\notin M_e\).
\end{lemma}

\begin{proof}
Put $N_G(v)=\{x,y,z\}$ and $H=G^v$.
A join rooted at $v_x$ contains $v_xx,v_xv_y,v_xv_z$.
Hence it contains neither $v_yy$ nor $v_zz$ nor $v_yv_z$.
Deleting its two triangle edges and contracting the new triangle
therefore gives a perfect matching of $G$ containing $vx$.
Conversely, a perfect matching $M$ of $G$ containing $vx$ gives
such a join with edge set
\(
(M\setminus\{vx\})\cup\{v_xx,v_xv_y,v_xv_z\}.
\)
By \cref{prop:loss-formula}, every lonely edge of $H$ is the transfer
of a lonely edge of $G$, since $v$ is zero-loss expandable.
A transferred edge from an old lonely edge incident with $v$ either
is incident with $v_x$ or equals $v_yv_z$, which the join cannot contain.
Thus a $v_x$-rooted join contains a lonely edge not incident with $v_x$
if and only if there is an $e\in L(G)$ not incident with $v$ such that
$vx\in M_e$. Taking the negation gives the assertion.
\begin{comment}
    Note that a join rooted at $v_x$ must contain all three edges incident with $v_x$, namely $v_xx,v_xv_y,v_xv_z$. Hence the other two vertices of the new triangle are already matched, and the remaining part corresponds exactly to a perfect matching of $G$ that uses $vx$. If a lonely edge $e$ is not incident with $v$, then such a join contains $e$ in $G^v$ if and only if the unique perfect matching $M_e$ of $G$ containing $e$ uses $vx$.

If a lonely edge is incident with $v$, its transferred edge lies inside the new triangle: either it is incident with the root $v_x$ and hence does not count toward the loss set, or it is the opposite edge $v_yv_z$ of $vx$, which a $v_x$-rooted join cannot use. Thus only lonely edges not incident with $v$ need to be checked, proving the claim.
\end{comment}
\end{proof}

\subsection{R-join and its property}
\begin{comment}
 To prove that no additional expandable vertices suddenly appear after one of the three local operations, we need the following strengthened notion used only in the proof. It is not part of the actual construction.   
\end{comment}

\begin{definition}\label{def:R-join}
Let $R\subseteq V(G)$. An $R$-join is a spanning subgraph $J$ of $G$ such that
\[
d_J(u)=3\quad(u\in R),\qquad d_J(u)=1\quad(u\notin R).
\]
Let
\(
\Act(G)=\set{v\in V(G):\Sigma_G(v)=\varnothing}
\)
denote the set of all zero-loss expansion vertices of $G\in\U_2$. For $A\subseteq V(G)$, we say that $A$ has the \emph{noninterference property} if, for every nonempty $R\subseteq A$ and every lonely edge $e$ with $V(e)\cap R=\varnothing$, there is no $R$-join containing $e$ (when $R=\set{v}$, this says exactly that $v$ is a zero-loss expansion vertex).
\end{definition}
%Allowing arbitrary nonempty $R$ serves only to ensure that after expanding one marked vertex, the other marked vertices remain expandable.

\begin{proposition}\label{prop:mark-propagation}
Suppose $G\in\U_2$ and suppose $\Act(G)$ has the noninterference property. Choose $v\in \Act(G)$ and put $H=G^v$. Define
\(
B_v=\set{u\in \Act(G)\setminus\set{v}:uv\in L(G)}
\)
and
\(
T_v=\set{v_x:x\in N_G(v),\;vx\notin M_e
\text{ for every }e\in L(G)\text{ not incident with }v}.
\)
Then \(\Act(H)=(\Act(G)\setminus(\set{v}\cup B_v))\cup T_v,\)
and $\Act(H)$ again has the noninterference property.
\end{proposition}

\begin{proof}
First determine which old vertices remain zero-loss expandable in $G^v$. If $u\notin \Act(G)$, choose a $u$-rooted join containing a lonely edge not incident with $u$. At $v$ this join uses one incident edge; replacing that edge by the corresponding external edge and opposite edge gives the same witness in $H$, so $u$ remains not zero-loss expandable.

Now let $u\in \Act(G)\setminus\set{v}$. If $uv$ is a lonely edge, a $u$-rooted join exists by \cite[Lemma~3.4]{GMMRW2026}, and every such join contains $uv$. Lifting it to $H$ makes the join contain the transferred opposite edge of $uv$, which is no longer incident with $u$, so $u$ ceases to be zero-loss expandable. These are exactly the vertices in $B_v$. Suppose instead that $uv$ is not lonely. We claim that every $u$-rooted join in $H$ contains no lonely edges not incident with $u$. %dose not contain any lonely edge incident with $u$.  

Suppose, to the contrary, that $H$ has a $u$-rooted join containing a lonely edge not incident with $u$. On the new triangle, the join has one of two forms. If it uses one triangle edge, contraction yields a $u$-rooted join in $G$ containing the corresponding lonely edge. This edge is still not incident with $u$, because $uv$ is not lonely, contradicting $u\in\Act(G)$. If it uses no triangle edge, then all three external edges are used, and contraction yields a $\set{u,v}$-join in $G$ containing the same lonely edge. That lonely edge is incident with neither $u$ (by its choice) nor $v$ (by \cref{cor:new} and the absence of selected triangle edges); this contradicts noninterference of $ \Act(G)$. 

By the claim above, a vertex \(u\) such that \(uv\) is not a lonely edge remains zero-loss expandable in \(G^v\). By \cref{lem:new-triangle-active}, the new-triangle vertices in \(T_v\) belong to \(\Act(H)\). This proves that
\(
\Act(H)=(\Act(G)\setminus(\set{v}\cup B_v))\cup T_v.
\)

It remains to prove noninterference of $\Act(H)$. Suppose, to the contrary, that for some nonempty $R'\subseteq \Act(H)$ there is an $R'$-join containing a lonely edge not incident with $R'$. Let $T$ be the new triangle, and let $k$ be the number of roots of $R'$ in $T$.
If $k=2$, both roots must select all three of their incident edges, so the third new vertex receives at least two selected triangle edges and has degree at least $2$, impossible for a non-root.

If $k=0$, then on the new triangle the join either uses one triangle edge, in which case contraction gives an $R'$-join in $G$, or it uses all three external edges, in which case contraction gives an $(R'\cup\set{v})$-join. In the first case, if the contained lonely edge is the transfer of $vx\in L(G)$, then $x\notin R'$: otherwise $x\in B_v$, whereas $B_v\cap\Act(H)=\varnothing$. Thus the corresponding old lonely edge is still nonincident with $R'$. In the second case, the contained lonely edge lies outside the new triangle and, by \cref{cor:new}, is not an external edge, so it remains nonincident with $R'\cup\{v\}$. Both cases contradict noninterference of $\Act(G)$.

If $k=1$, let the unique new root be $v_x$. The contained lonely edge lies outside the new triangle: its only edge not incident with $v_x$ cannot be selected in a $v_x$-rooted join, and external edges are not lonely. Delete its two triangle edges and contract the new triangle. This yields an $(R'\setminus\set{v_x})$-join in $G$ that is forced to use $vx$. If the old root set is nonempty, this again contradicts noninterference of $\Act(G)$. If the old root set is empty, we obtain a perfect matching containing the lonely edge and using $vx$, which must be the corresponding unique perfect matching $M_e$; this contradicts $v_x\in T_v$.

If $k=3$, the three new roots force all three triangle edges and all three external edges to be selected. Contracting the new triangle gives a $((R'\setminus V(T))\cup\set{v})$-join in $G$. Since the contained lonely edge is incident with none of the three new roots, it cannot be the transfer of a lonely edge incident with $v$, and after contraction it remains nonincident with every root. This again contradicts noninterference of $\Act(G)$. All four cases are impossible, so $\Act(H)$ remains noninterfering.
\end{proof}

\subsection{Two standard markings and three expansion operations}\label{sec:three-operations}

\begin{definition}\label{def:equal-standard-mark}
Let $G\in\mathcal U_2$. A marked object
\(
\mathfrak G=(G;ab,de;abc,def)
\)
is called an \emph{equal-matching standard marking} if the following hold:
\begin{enumerate}[label=\textup{(E\arabic*)}]
\item $L(G)=\set{ab,de}$ and $M_{ab}=M_{de}$;%, and the common unique perfect matching of the two lonely edges is denoted by $M$;
\item $abc$ and $def$ are triangles of $G$;
\item
\(\Act(G)=\set{a,b,d,e}.\)
\end{enumerate}
\begin{comment}{The order of the two lonely edges, the order of the two endpoints of each edge, and the third triangle vertices $c,f$ are all part of the marking data. Interchanging these orders merely gives another standard marking of the same unlabeled graph.}
\end{comment}
\end{definition}

\begin{lemma}\label{lem:equal-mark-noninterference}
If $\mathfrak G=(G;ab,de;abc,def)$ is an equal-matching standard marking, then
$\Act(G)=\set{a,b,d,e}$ has the noninterference property.
\end{lemma}

\begin{proof}
Let $R$ be an arbitrary nonempty subset of $\Act(G)$. Consider an $R$-join of $G$. Note that $abc$ and $def$ are triangles and $\Act(G)=\set{a,b,d,e}$. By the definition of an $R$-join, we have $|\{a, b, c\} \cap R| \leq 1$ and $|\{d, e, f\} \cap R| \leq 1$. If $|\{a, b, c\} \cap R| = 1$ and $|\{d, e, f\} \cap R| = 1$, then no lonely edge is nonincident to $R$. This implies that $\Act(G)$ satisfies the noninterference property. Otherwise, $|R \cap \{a, b, d, e\}| = 1$. Without loss of generality, assume $R = \{a\}$. Then an $R$-join containing $de$ is exactly a rooted join at the unique root that includes a nonincident lonely edge, which contradicts the fact that the root is in $\Act(G)$.
\begin{comment}
 and let the lonely edge $e_0$ be nonincident with $R$. If $e_0=de$, then $R\subseteq\set{a,b}$; if $e_0=ab$, then $R\subseteq\set{d,e}$. The two cases are symmetric, so consider the first.
If $|R|=1$, then an $R$-join containing $de$ is precisely a rooted join at the unique root that contains a nonincident lonely edge, contradicting that the root belongs to $\Act(G)$. If $R=\set{a,b}$, then because both $a$ and $b$ are roots, the $R$-join must use both $ac$ and $bc$. Since $c\notin R$, the degree of $c$ in this spanning subgraph is then at least $2$, contradicting the requirement that every non-root have degree $1$. Hence no such $R$-join exists.   
\end{comment}
\end{proof}

Let $\mathfrak G=(G;ab,de;abc,def)$ be an equal-matching standard marking, and consider the endpoint $a$ of the lonely edge $ab$. Let $p$ be the third neighbor of $a$ outside the triangle $abc$. 
\begin{comment}
    Expand at $a$ and label the new triangle by
\[
\Delta_a=a_ba_ca_p,
\qquad a_b\sim b,\quad a_c\sim c,\quad a_p\sim p.\]
\end{comment}

\begin{definition}\label{def:operation-E}
In the above situation, let $H=G^a$. Update the old lonely edge and the marking by
\[
ab\longmapsto a_ca_p,
\qquad
\mathfrak G^a=(H;a_ca_p,de; a_ca_pa_b,def).
\]
This is called \emph{endpoint advance} at $a$, or Operation E at $a$. The operations at endpoints $b,d,e$ are defined by the completely symmetric rule.
\end{definition}

\begin{figure}[htbp]
\centering
\includegraphics[width=0.55\textwidth]{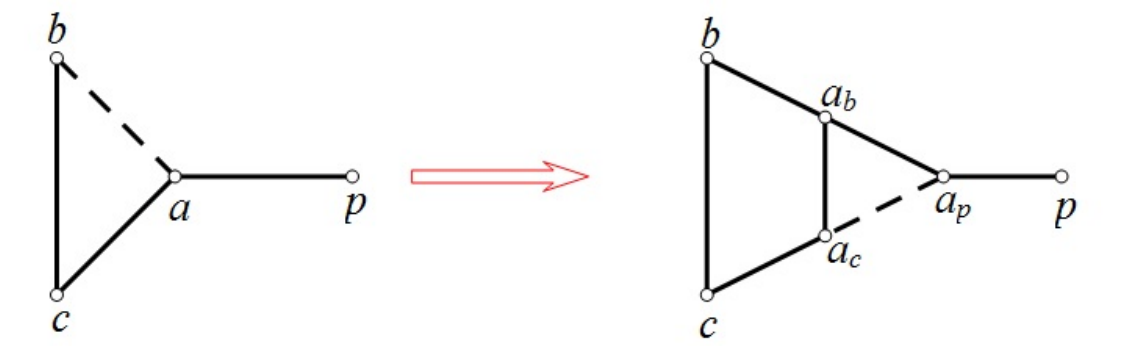}
\caption{Operation E. %The dashed lonely edge $ab$ is advanced to the opposite edge $a_ca_p$ in the new triangle.
}
\label{fig:operation-E}
\end{figure}

\begin{lemma}\label{lem:operation-E}
If $\mathfrak G=(G;ab,de;abc,def)$ is an equal-matching standard marking, then after applying Operation E at any endpoint of a lonely edge, the resulting $\mathfrak G^a$ is again an equal-matching standard marking. More precisely,
\(
L(H)=\set{a_ca_p,de},
\Act(H)=\set{a_c,a_p,d,e}.
\)
\end{lemma}

\begin{proof}
By Definition~\ref{def:equal-standard-mark}, $a\in\Act(G)$. This means that  $\Sigma_G(a)=\varnothing$. By \cref{prop:loss-formula}, $H=G^a\in\U_2$ and
\(
L(H)=\set{\lambda_a(ab),\lambda_a(de)}
     =\set{a_ca_p,de}.
\)
Recall that $M_{ab}=M_{de}$. By Lemma~\ref{lem:matching-type-preserved}, $M_{a_ca_p}=M_{de}$.  

We now compute the zero-loss vertex set of $H$ exactly. Since
\(
\Act(G)=\set{a,b,d,e}
\)
and $ab\in L(G)$, the set appearing in Proposition~\ref{prop:mark-propagation} is
\( B_a=\set{b}.\)
The zero-loss vertices in the new triangle are given by $T_a$. The only lonely edge not incident with $a$ is $de$, and $M_{de}=M_{ab}$. 
Clearly, $ac, ap\notin M_{de}$. 
By the definition of $T_a$,
\(T_a=\set{a_c,a_p}.\)
By Lemma~\ref{lem:equal-mark-noninterference}, $\Act(G)$ is noninterfering, so Proposition~\ref{prop:mark-propagation} applies and yields
\(
\Act(H)=(\Act(G)\setminus(\set{a}\cup B_a))\cup T_a=\set{a_c,a_p,d,e}.
\)
The new lonely edge $a_ca_p$ lies in the new triangle $a_ba_ca_p$, while the other lonely edge $de$ and triangle $def$ are unchanged. Hence all three conditions of Definition~\ref{def:equal-standard-mark} hold, so $\mathfrak G^a$ is an equal-matching standard marking.
\end{proof}

\begin{definition}\label{def:unequal-standard-mark}
A marked object
\(
\mathfrak G=(G;x,v,w,y;\mathcal S,\rho)
\)
is called an \emph{unequal-matching standard marking} if the following conditions hold.

\begin{enumerate}[label=\textup{(D\arabic*)}]
\item $G\in\U_2^{\ne}$, and
\(
x-v-w-y
\)
is a path with  
\(
L(G)=\{xv, wy\}.
\)
\item
\(
\Act(G)=\set{v,w}\cup\mathcal S, \mathcal S\subseteq V(G)\setminus\{x,v,w,y\},
\)
and $\Act(G)$ has the noninterference property.
\begin{comment}
 \item If
\(
N_G(v)=\set{x,w,t}, N_G(w)=\set{v,y,s},
\)
then the two unique perfect matchings satisfy
\begin{equation}\label{eq:center-matching-pattern}
\begin{aligned}
M_{xv}\cap\delta(v)&=\set{xv},&
M_\beta\cap\delta(v)&=\set{vt},\\
M_\beta\cap\delta(w)&=\set{wy},&
M_{xv}cap\delta(w)&=\set{ws}.
\end{aligned}
\end{equation}
In particular, the center edge $vw\notin M_{xv}\cup M_\beta$.   
\end{comment}

\item $\rho:\mathcal{S}\to V(G)$ is a mapping satisfying the following. For each $u\in\mathcal{S}$, put $p=\rho(u)$. Then $p\in N_G(u)\setminus\Act(G)$ and, writing $N_G(u)=\{p,q,r\}$, we have
%\begin{equation}%\label{eq:tip-matching-pattern}
\(\left\{M_{xv}\cap\delta(u),\,M_{wy}\cap\delta(u)\right\}
=\left\{\{uq\},\{ur\}\right\},\)
%\end{equation}
where $\delta(u)$ denotes the set of edges incident to vertex $u$.
\end{enumerate}
\end{definition}

The vertices $v,w$ are called the center vertices, and $vw$ the center edge in $\mathfrak G$. For (D3), we write $u\to p$, call $u$ a tip, and call $p$ its direction. The requirement $\rho(u)\notin\Act(G)$ explicitly ensures that an arrow never points to another zero-loss vertex. This eliminates the ambiguity that expanding one tip might simultaneously delete the direction target of another tip.

Let $\mathfrak G=(G;x,v,w,y;\mathcal S,\rho)$ be an unequal-matching standard marking, with $N_G(v)=\set{x,w,t}$. Expand at the center vertex $v$, labeling the new triangle $v_x,v_w,v_t$. The transferred opposite edge of the old lonely edge $xv$ is $v_wv_t$.

\begin{definition}\label{def:operation-C}
Let $H=G^v$. Define the new center path, tip set, and direction map by
\[
v_t-v_w-w-y,\qquad
\mathcal S'=\mathcal S\cup\set{v_x},
\qquad
\rho'(v_x)=v_t,
\qquad
\rho'|_{\mathcal S}=\rho.
\]
This gives the marking
\(
\mathfrak G^v=(H;v_t,v_w,w,y;\mathcal S',\rho').
\)
The operation is called \emph{center expansion} at $v$, or Operation C at $v$. The rule at the right center $w$ is completely symmetric.
\end{definition}

\begin{figure}[htbp]
\centering
\includegraphics[width=0.65\textwidth]{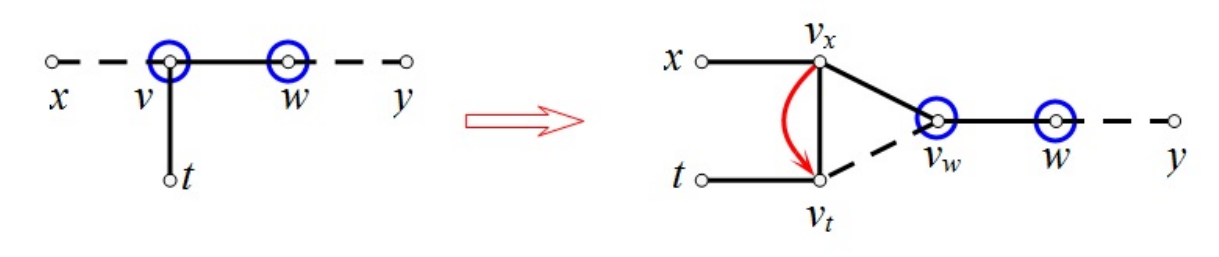}
\caption{Operation C: the new centers are $v_w,w$, and a new tip $v_x\to v_t$ is created in $H$.}
\label{fig:operation-C}
\end{figure}

\begin{lemma}\label{lem:operation-C}
Applying Operation C at either center vertex of an unequal-matching standard marking again yields an unequal-matching standard marking and 
\(
\Act(H)=\set{v_w,w}\cup\mathcal S\cup\set{v_x}.
\)
Thus all old tips are retained, and exactly one new tip, $v_x\to v_t$, is added.
\end{lemma}

\begin{proof}
Note that $v\in\Act(G)$. By \cref{prop:loss-formula}, \(L(H)=\set{v_wv_t,wy},\)
and Lemma~\ref{lem:matching-type-preserved} gives $H\in\U_2^{\ne}$. Since $x\notin\Act(G)$, it follows from the definition of $B_v$ that $B_v=\varnothing$.
Note that the unique lonely edge not incident with $v$ is $wy$. Because $M_{xv}\not= M_{wy}$, $xv\notin M_{wy}$. Clearly, $vw\notin M_{wy}$. Therefore, the set \(T_v\)
as in Proposition~\ref{prop:mark-propagation} is \(\set{v_x,v_w}\). It follows that
\(
\Act(H)
=(\Act(G)\setminus\set{v})\cup\set{v_x,v_w}
=\set{v_w,w}\cup\mathcal S\cup\set{v_x},
\)
 and this set remains noninterfering. To complete the proof, it suffices to show $\mathfrak G^v$ satisfies the third  condition in \cref{def:unequal-standard-mark}.

\begin{comment} Next verify the center-matching pattern. 
    Since $M_\alpha$ uses $vx$ at $v$, its natural lift uses $v_xx$ and $v_wv_t$ in the new triangle. Since $M_\beta$ uses $vt$ at $v$, its natural lift uses $v_tt$ and $v_xv_w$. Thus at the new left center $v_w$, the matching $M_{v_wv_t}$ uses the new lonely edge $v_wv_t$, whereas $M_{wy}$ uses the third edge $v_wv_x$; the center edge $v_ww$ belongs to neither matching. At the old right center $w$, the original matching relations are unchanged except that the old edge $wv$ is replaced by the external edge $wv_w$, so~\eqref{eq:center-matching-pattern} still holds.
\end{comment}
At the new vertex $v_x$, the two unique matchings use $v_xx$ and $v_xv_w$, respectively, while $v_xv_t$ belongs to neither; moreover $v_t\notin\Act(H)$. Thus $v_x\to v_t$ satisfies the conditions. Finally consider an old tip $u\to\rho(u)$. Since $\rho(u)\notin\Act(G)$ but $v\in\Act(G)$, no old arrow can point to $v$. If $u$ is not adjacent to $v$, its local matching pattern is unchanged; if $u$ is adjacent to $v$, the old edge $uv$ is simply replaced by the corresponding external edge $uv_u$, and the relation describing which nondirectional edge is used at $u$ by each matching is unchanged. Proposition~\ref{prop:mark-propagation} also shows that every old non-zero-loss vertex remains non-zero-loss, so every old arrow target remains outside the zero-loss vertex set. 
\end{proof}

Let $\mathfrak G=(G;x,v,w,y;\mathcal S,\rho)$ be an unequal-matching standard marking. Choose a tip $u\in\mathcal S$, put $p=\rho(u)$, and write
\(
N_G(u)=\set{p,q,r}.
\)
After exchanging the names of $q,r$ if necessary, assume that $M_{xv}$ uses $uq$ at $u$ and $M_{wy}$ uses $ur$. Expand at $u$ to the new triangle $u_pu_qu_r$.

\begin{definition}\label{def:operation-T}
Let $H=G^u$ and put
\(\mathcal S'=(\mathcal S\setminus\set{u})\cup\set{u_p}.\) Define
\(\rho'(u_p)=p, \rho'|_{\mathcal S\setminus\set{u}}=\rho.
\)
The center path remains $x-v-w-y$. The resulting marking
\(
\mathfrak G^u=(H;x,v,w,y;\mathcal S',\rho')
\)
is called \emph{tip extension}, or Operation T.
\end{definition}

\begin{figure}[htbp]
\centering
\includegraphics[width=0.6\textwidth]{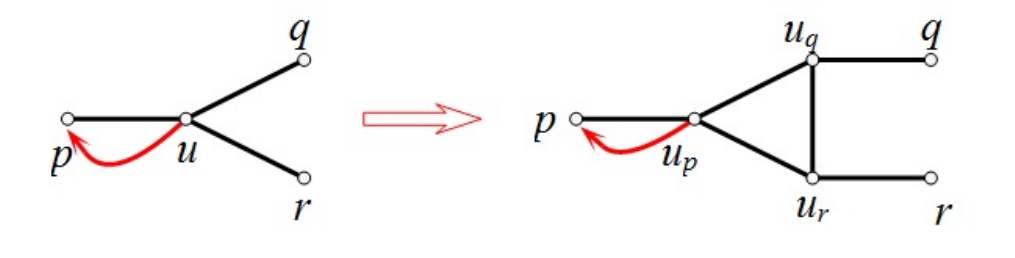}
\caption{Operation T. The old tip $u\to p$ is replaced by the new tip $u_p\to p$.}
\label{fig:operation-T}
\end{figure}

\begin{lemma}\label{lem:operation-T}
Applying Operation T at any tip of an unequal-matching standard marking again yields an unequal-matching standard marking, and
\(
\Act(H)=\set{v,w}\cup(\mathcal S\setminus\set{u})\cup\set{u_p}.
\)
The center vertices and all other tips remain unchanged, and the number of expandable vertices is unchanged.
\end{lemma}

\begin{proof}
By \cref{def:unequal-standard-mark}, the tip $u$ does not belong to the four vertices of the center path, so neither lonely edge is incident with $u$, and $B_u=\varnothing$. By the tip condition, $up$ belongs to neither $M_{xv}$ nor $M_{wy}$, while $uq$ and $ur$ are used by the two matchings, respectively. Hence by Lemma~\ref{lem:new-triangle-active}, $u_p$ is zero-loss expandable in the new triangle, that is,
\(
T_u=\set{u_p}.
\)
Proposition~\ref{prop:mark-propagation} gives
\(
\Act(H)=(\Act(G)\setminus\set{u})\cup\set{u_p}
       =\set{v,w}\cup(\mathcal S\setminus\set{u})\cup\set{u_p},
\)
and guarantees that this set remains noninterfering. By \cref{prop:loss-formula} and Lemma~\ref{lem:matching-type-preserved}, $H\in\U_2^{\ne}$, and the two lonely edges and the center path are unchanged.

After natural lifting, $M_{xv}$ uses $u_qq$ and $u_pu_r$ in the new triangle, while $M_{wy}$ uses $u_rr$ and $u_pu_q$. Thus at $u_p$ the two matchings use two different triangle edges, and the external edge $u_pp$ belongs to neither; the new arrow is therefore $u_p\to p$. Since $p\notin\Act(G)$ originally, the formula above shows that $p\notin\Act(H)$ as well.

For any other old tip $z\to\rho(z)$, arrow targets lie outside $\Act(G)$, so $\rho(z)\ne u$. Hence expanding $u$ does not delete the target of any other arrow. If $z$ is adjacent to $u$, the corresponding nondirectional edge is merely replaced by an external edge and the matching-use pattern is preserved; otherwise the local structure is unchanged. The same argument applies at the center vertices: if $u$ is the third neighbor of a center, the relevant old matching edge is simply rerouted to the corresponding new vertex, so the local matching pattern at each center is preserved. Hence all standard-marking conditions are preserved.
\end{proof}

\begin{corollary}\label{cor:branching-law}
In an equal-matching standard marking, the zero-loss vertices are always exactly the four endpoints of the two lonely edges, so every step has exactly four possible Operation E choices.

In an unequal-matching standard marking with $k$ tips, the zero-loss vertices are exactly the two center vertices together with those $k$ tips. Operation C at either center increases the number of tips by $1$, whereas the other $k$ Operation T choices merely extend existing tips and leave the number of tips unchanged.
\end{corollary}

\begin{remark}\label{rem:no-circularity}
Up to this point we have proved only the following conditional statement: \emph{once} a graph carries an equal-matching or unequal-matching standard marking, the allowed local operations produce another standard marking of the same type. We have not yet proved that the lonely edges of every $G\in\U_2^{=}$ lie in triangles, nor that every $G\in\U_2^{\ne}$ has a center path and tips. These general conclusions will be derived in Section~\ref{sec:main-proof} from the primitive-kernel base cases and kernel-reduction induction.
\end{remark}

\subsection{Standard-marking base cases for the primitive kernels}\label{sec:kernel-marking}

In this subsection we do not assume that a general graph in $\U_2$ already has the standard structure. Instead, we verify one by one that every primitive kernel admits at least one standard marking in the sense of Subsection~\ref{sec:three-operations}. These kernels will serve as the induction bases in Subsection~\ref{sec:main-proof}. In the figures, dashed edges are lonely edges, double circles indicate center vertices, the tail of an arrow is a tip, and the head of the arrow is its direction.

Using the integer labels from Subsection~\ref{sec:finite-branch}, all three kernels are of unequal-matching type. Their center paths and tips are
\begin{align*}
D_0:&\quad 10-1-8-4,
&0\to5,\quad 2\to3,\quad 9\to4,\quad 11\to10;\\
D_1:&\quad 3-1-8-4,
&2\to3,\quad 9\to4,\quad 13\to5;\\
D_2:&\quad 1-3-6-7,
&2\to1,\quad 9\to5,\quad 12\to7.
\end{align*}
In each row, the first and last edges of the four-vertex path are the lonely edges and the two internal vertices are the centers.

 \begin{figure}[htbp]
\centering
\includegraphics[width=0.88\textwidth]{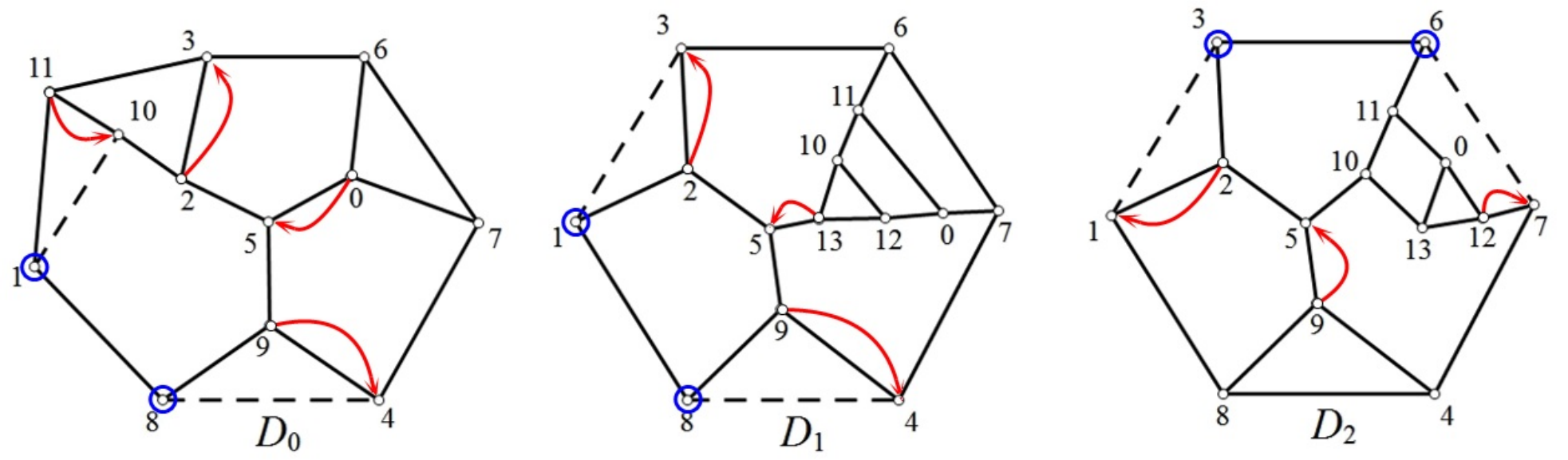}
\caption{One unequal-matching standard marking for three non-extended-prism kernels. %Ordinary graph edges are undirected; the additional arrows indicate expansion directions only.
}
\label{fig:D-marks}
\end{figure}

Use the natural labels $U=u_1u_2u_3$ and $X^j=x_1^jx_2^jx_3^j$. Standard markings for the three equal-matching kernels are listed below. The column ``containing triangles'' explicitly records the structural data from Definition~\ref{def:equal-standard-mark} that cannot be omitted in advance.

\begin{center}
\begin{tabular}{c c c c}
\toprule
Kernel & Two lonely edges & Containing triangles & $\Act(G)$\\
\midrule
$P(1122)$ & $u_1u_2,\quad x_1^4x_2^4$ & $U,\quad X^4$ & $u_1,u_2,x_1^4,x_2^4$\\
$P(1221)$ & $u_1u_3,\quad x_1^4x_3^4$ & $U,\quad X^4$ & $u_1,u_3,x_1^4,x_3^4$\\
$P(1231)$ & $u_1u_3,\quad x_1^4x_2^4$ & $U,\quad X^4$ & $u_1,u_3,x_1^4,x_2^4$\\
\bottomrule
\end{tabular}
\end{center}

Standard markings for the four unequal-matching kernels are as follows.
\begin{center}
\begin{tabularx}{\textwidth}{c X X}
\toprule
Kernel & Three-edge center path & Tips and directions\\
\midrule
$P(1121)$ & $x_1^2-x_3^2-x_3^3-x_2^3$
& $u_1\to x_1^2,\quad x_2^4\to x_2^3$\\
$P(1212)$ & $x_1^1-x_3^1-x_3^2-x_2^2$
& $u_1\to x_1^1,\quad x_2^4\to x_2^2$\\
$P(1213)$ & $x_3^0-x_2^0-x_2^2-x_3^2$
& $u_3\to x_3^0,\quad x_3^1\to x_1^1,\quad x_3^4\to x_3^2$\\
$P(1223)$ & $x_3^0-x_2^0-x_2^3-x_1^3$
& $u_3\to x_3^0, x_3^1\to x_1^1, x_1^2\to x_3^2,x_1^4\to x_1^3$\\
\bottomrule
\end{tabularx}
\end{center}

\begin{lemma}\label{lem:finite-seed-certificate}
For the ten finite primitive kernels listed in Figure~\ref{fig:D-marks} and the two tables above, the following all hold:
\begin{enumerate}[label=\textup{(\alph*)}]
\item the listed lonely edges and matching types are correct;
\item the listed vertices are exactly all zero-loss expansion vertices;
\item the three equal-matching graphs satisfy Definition~\ref{def:equal-standard-mark};
\item the seven unequal-matching graphs satisfy Definition~\ref{def:unequal-standard-mark}; in particular, the matching direction at every tip, inactivity of every arrow target, and noninterference of the zero-loss vertex set all hold.
\end{enumerate}
\end{lemma}

\begin{proof}
For each of the ten graphs, we enumerate its perfect matchings and
count the occurrences of each edge to determine its lonely edges
and their unique perfect matchings. Comparing these two matchings
determines the matching type. We compute the complete zero-loss
vertex set using \cref{lem:deletion-criterion}, and check the listed
triangles, center paths, and tip directions directly.
For each nonempty subset $R$ of the proposed zero-loss vertex set
and each lonely edge $e$ nonincident with $R$, we select all edges
incident with $R$ together with $e$. If a non-root has degree greater
than $1$, no such $R$-join exists; otherwise, the remaining uncovered
vertices must have a perfect matching. Exhaustive matching enumeration
excludes this possibility in every case. These checks establish all
four assertions; the computations and certificates are provided
in~\cite{qianli}.
\end{proof}

For the next two lemmas, let $G=P(12^rt)$, where $r\ge2$ and
$t\in\{1,3\}$, and put $q=4-t$.
To specify the alternating edge sets, including their boundary edges,
write
\[e_0^1=u_1x_1^1,\  e_0^3=x_3^0x_3^1, \ e_r^t=x_t^rx_t^{r+2}, \ e_r^q=x_q^rx_q^{r+1}, \ e_i^j=x_j^ix_j^{i+1} (1\le i<r,\ j\in\{1,3\}).\] 
The pair $\{e_i^1,e_i^3\}$ is the pair of rail edges at cut $i$.
A rail-edge selection is called \emph{alternating} if it selects one
edge at each consecutive cut and switches rails at each step.
Put $\eta_i=1$ for even $i$ and $\eta_i=3$ for odd $i$, and define
\[
\mathcal A_{a,b}=\{e_i^{\eta_i}:a\le i\le b\},\qquad
\mathcal B_{a,b}=\{e_i^{4-\eta_i}:a\le i\le b\},\qquad
\mathcal R=\{x_1^kx_3^k:1\le k\le r\}.
\]
Here $0\le a\le b\le r$, and an interval with $a>b$ denotes the
empty set. Thus $\mathcal R$ contains only the internal rungs.
If neither vertex of an internal rung is a root and exactly one
rail edge is selected at one adjacent cut, the degree requirements
exclude the rung and force the other rail at the next cut.
Consequently, alternating selections propagate through every
segment containing no roots.

Let %Define the alternating candidate vertices and their opposite vertices in the same ladder cells by
% Match the numbering in the screenshots: (16), followed by (17)--(20).

\begin{equation}\label{eq:alternating-candidates}
s_k=
\begin{cases}
x_3^k,&k\text{ is odd},\\
x_1^k,&k\text{ is even},
\end{cases}
\qquad
\bar s_k=
\begin{cases}
x_1^k,&k\text{ is odd},\\
x_3^k,&k\text{ is even}.
\end{cases}
\end{equation}
Here $1\le k\le r$. Clearly,
$\{s_k,\bar s_k\}=\{x_1^k,x_3^k\}$.

\begin{lemma}
\label{lem:infinite-common-mark}
Let $r\ge2$ and $t\in\{1,3\}$. If $t=1$ and $r$ is even, or $t=3$ and $r$ is odd, then
\[
L(P(12^rt))=\set{u_1u_3,x_1^{r+2}x_3^{r+2}},\ 
M_{u_1u_3}=M_{x_1^{r+2}x_3^{r+2}},\ 
\Act(P(12^rt))=\set{u_1,u_3,x_1^{r+2},x_3^{r+2}}.
\]
\end{lemma}

\begin{proof}
By \cref{lem:ladder-lonely} and Corollary~\ref{cor:infinite-kernel-type},
we have $L(P(12^rt))=\set{u_1u_3,x_1^{r+2}x_3^{r+2}}$ and
$M_{u_1u_3}=M_{x_1^{r+2}x_3^{r+2}}$.
To complete the proof, it suffices to show that
\[
\Act(P(12^rt))=\set{u_1,u_3,x_1^{r+2},x_3^{r+2}}.
\]

First, we show that $\{u_1,u_3,x_1^{r+2},x_3^{r+2}\}\subseteq\Act(G)$
by contradiction. Suppose that there exists a $y$-rooted join $J_y$
containing a lonely edge not incident with $y$, where
$y\in\{u_1,u_3,x_1^{r+2},x_3^{r+2}\}$.

Since $u_3$ is the root of $J_{u_3}$,
$\{u_1u_3,u_2u_3,u_3x_3^0\}\subseteq E(J_{u_3})$.
The degree requirement at $x_2^0$ then forces
$x_2^0x_2^{r+1}\in E(J_{u_3})$.
On the other hand, the nonincident lonely edge in $J_{u_3}$ must be
$x_1^{r+2}x_3^{r+2}$, which forces
$x_2^{r+1}x_2^{r+2}\in E(J_{u_3})$.
Thus the non-root $x_2^{r+1}$ has degree at least $2$, a contradiction.

For $J_{u_1}$, the root forces
$u_1u_2,u_1u_3,u_1x_1^1\in E(J_{u_1})$.
The nonincident lonely edge $x_1^{r+2}x_3^{r+2}$ forces
$x_2^{r+1}x_2^{r+2}$ and excludes $x_2^0x_2^{r+1}$.
Hence $x_2^0x_3^0$ is selected, and $x_3^0x_3^1$ is excluded.
At the right end, the non-root $x_q^{r+1}$ can therefore only be
covered by $e_r^q$.
However, the entrance edge $e_0^1=u_1x_1^1$ forces the alternating
selection $\mathcal A_{0,r}$ through the root-free internal rungs.
In the present parity cases $\eta_r=t$, so the selected exit is
$e_r^t$, not $e_r^q$. This is a contradiction.

Now let the root be $x_q^{r+2}$ and suppose that the join contains
$u_1u_3$. The left end forces $u_2x_2^0,x_3^0x_3^1$ and excludes
$x_2^0x_2^{r+1}$. The three edges at the root cover
$x_q^{r+1},x_2^{r+2},x_t^{r+2}$. Consequently, none of the three
neighbors of $x_2^{r+1}$ can be joined to it, a contradiction.
Finally, if the root is $x_t^{r+2}$, the same left-end forcing
excludes $x_2^0x_2^{r+1}$ and gives entrance $e_0^3$.
The root forces $e_r^t$ and its two triangle edges, after which
$x_2^{r+1}x_q^{r+1}$ is forced and $e_r^q$ is excluded.
But the entrance $e_0^3$ propagates as $\mathcal B_{0,r}$ and,
since $\eta_r=t$, exits on rail $q$, again a contradiction.
Thus all four claimed vertices are zero-loss expandable.

To prove that no other vertex is zero-loss expandable, we construct
a $y$-rooted join containing a lonely edge not incident with $y$
for every
$y\in V(G)\setminus\{u_1,u_3,x_1^{r+2},x_3^{r+2}\}$.
In each row below, the edge set of the join is the union of
$\delta(y)$, the listed lonely edge, and the edge sets in the last
two columns. Repeated edges are included only once, and $1\le k\le r$.

\begin{center}
\scriptsize
\begin{tabular}{c c l l}
\toprule
Root $y$ & Lonely edge & Other boundary edges & Internal rungs and rail edges\\
\midrule
$x_2^0$ & $u_1u_3$ & $x_q^{r+1}x_q^{r+2},\ x_2^{r+2}x_t^{r+2}$
& $\mathcal R$\\[1ex]
$x_q^{r+1}$ & $u_1u_3$ & $u_2x_2^0,\ x_2^{r+2}x_t^{r+2}$
& $\mathcal B_{0,r}$\\[1ex]
$x_2^{r+2}$ & $u_1u_3$ & $u_2x_2^0$
& $\mathcal B_{0,r}$\\[1ex]
$s_k$ & $u_1u_3$ & $u_2x_2^0,\ x_q^{r+1}x_q^{r+2},\ x_2^{r+1}x_2^{r+2}$
& $\mathcal B_{0,k-1}\cup\mathcal A_{k,r}$\\[1ex]
$\bar s_k$ & $x_1^{r+2}x_3^{r+2}$ & $u_1x_1^1,\ u_2x_2^0,\ u_3x_3^0,\ x_2^{r+1}x_2^{r+2}$
& $\mathcal A_{0,k-1}\cup\mathcal B_{k,r}$\\[1ex]
$u_2$ & $x_1^{r+2}x_3^{r+2}$ & $x_2^{r+1}x_2^{r+2}$
& $\mathcal B_{0,r}$\\[1ex]
$x_3^0$ & $x_1^{r+2}x_3^{r+2}$ & $u_1u_2,\ x_2^{r+1}x_2^{r+2}$
& $\mathcal B_{0,r}$\\[1ex]
$x_2^{r+1}$ & $x_1^{r+2}x_3^{r+2}$ & $u_1u_2,\ u_3x_3^0$
& $\mathcal R$\\
\bottomrule
\end{tabular}
\end{center}

For an internal root $s_k$ or $\bar s_k$, the selections on cuts
$k-1$ and $k$ lie on the root's rail, and its rung is supplied by
$\delta(y)$. The root therefore has degree $3$, while the opposite
vertex has degree $1$. At every other internal rung, the two adjacent
selected rail edges lie on different rails, so both vertices have
degree $1$. For a row using $\mathcal R$, each internal vertex is
covered by its rung. The displayed boundary edges and $\delta(y)$
give degree $1$ at every remaining non-root vertex.
Thus every row defines the required join. The rows cover all
$2r+6$ vertices outside the claimed zero-loss set, completing the proof.
\end{proof}

\begin{lemma}
\label{lem:infinite-unequal-mark}
Let $r\ge2$ and $t\in\{1,3\}$. The extended prism $P(12^rt)$ admits an unequal-matching standard marking when $t=1$ and $r$ is odd, or $t=3$ and $r$ is even.
\end{lemma}
\begin{comment}
    If $t=1$ and $r$ is odd, or $t=3$ and $r$ is even, then
\(
L(P(12^rt))=\set{x_2^0x_3^0,x_2^{r+1}x_q^{r+1}}.
\)
The central three-edge path is
\begin{equation}\label{eq:infinite-center-path}
x_3^0-x_2^0-x_2^{r+1}-x_q^{r+1}.
\end{equation}
The complete collection of tips, with their directions, is
\begin{equation}\label{eq:infinite-tips}
u_3\to x_3^0,\qquad
s_k\to\bar s_k\ (1\le k\le r),\qquad x_q^{r+2}\to x_q^{r+1},
\end{equation}
and the complete zero-loss vertex set is
\begin{equation}\label{eq:unequal-infinite-active}
A=\set{x_2^0,x_2^{r+1},u_3,x_q^{r+2},s_1,\cdots,s_r}.
\end{equation}
The set $A$ is non-interfering. 
\end{comment}
\begin{proof}
We prove this result by showing that
$\mathfrak G=(G;x_3^0,x_2^0,x_2^{r+1},x_q^{r+1};\mathcal S,\rho)$
is an unequal-matching standard marking of $G=P(12^rt)$, where
$\mathcal S=\set{u_3,x_q^{r+2},s_1,\cdots,s_r}$ and
\[
\rho(u_3)=x_3^0,\qquad \rho(x_q^{r+2})=x_q^{r+1},\qquad
\rho(s_k)=\bar s_k\quad(1\le k\le r).
\]
The lonely edges and the distinct-matching type follow from
\cref{lem:ladder-lonely,cor:infinite-kernel-type}; explicitly,
\(
L(G)=\{x_2^0x_3^0,x_2^{r+1}x_q^{r+1}\}.
\)
Let $A=\{x_2^0,x_2^{r+1},u_3,x_q^{r+2},s_1,\ldots,s_r\}$.

\noindent\textbf{Claim:} $A$ is noninterfering.

\noindent\textit{Proof of Claim.}
Suppose, to the contrary, that there are a nonempty $R\subseteq A$,
an $R$-join $J$, and an edge $e\in L(G)\cap E(J)$ with
$V(e)\cap R=\varnothing$.
\begin{itemize}
\item Suppose $e=x_2^0x_3^0$.
Neither endpoint is a root, so $x_2^0x_2^{r+1}$,
$u_2x_2^0$, $u_3x_3^0$, and $x_3^0x_3^1$ are excluded from $J$.
In particular, $u_3,x_2^{r+1}\notin R$.
The degree requirements at $u_1,u_2,u_3$ then force
$u_2u_3,u_1x_1^1\in E(J)$.
If $x_q^{r+2}$ were a root, its three incident edges would cover
$x_q^{r+1},x_2^{r+2},x_t^{r+2}$. Since $x_2^0x_2^{r+1}$ is absent,
the non-root $x_2^{r+1}$ would have no available neighbor.
Thus $x_q^{r+2}\notin R$, and $R\subseteq\{s_1,\ldots,s_r\}$.

Let $s_j$ be the root with smallest index. Before cell $j$, all cells
contain no roots, so the entrance $e_0^1$ forces
$\mathcal A_{0,j-1}$. The selected edge at cut $j-1$ lies on rail
$\eta_{j-1}=4-\eta_j$, whereas the root $s_j=x_{\eta_j}^j$
requires the edge on rail $\eta_j$ at that cut. This is a contradiction.

\item Suppose $e=x_2^{r+1}x_q^{r+1}$.
Neither endpoint is a root. Hence $x_2^0x_2^{r+1}$,
$x_2^{r+1}x_2^{r+2}$, $e_r^q$, and $x_q^{r+1}x_q^{r+2}$
are excluded from $J$.
It follows that $x_2^0,x_q^{r+2}\notin R$.
If $u_3$ were a root, its three incident edges would cover
$u_1,u_2,x_3^0$, leaving $x_2^0$ with no available neighbor.
Thus $u_3\notin R$, and again $R\subseteq\{s_1,\ldots,s_r\}$.
The root-free right end forces $x_q^{r+2}x_2^{r+2}$ and $e_r^t$.

Let $s_j$ be the root with largest index. In the present parity cases,
$t=4-\eta_r$. Propagating backward from $e_r^t$ through the cells
with indices greater than $j$ selects the edge on rail $4-\eta_j$
at cut $j$. But the root $s_j=x_{\eta_j}^j$ requires the edge on
rail $\eta_j$ at that cut, a contradiction.
\end{itemize}
Both cases are impossible, proving the claim.

Taking singleton root sets in the noninterference property gives
$A\subseteq\Act(G)$. To prove equality, we construct the joins
listed below. As in the preceding proof, each join consists of
$\delta(y)$, the listed lonely edge, and the last two columns;
repeated edges are included only once, and $1\le k\le r$.

\begin{center}
\scriptsize
\begin{tabular}{c c l l}
\toprule
Root $y$ & Lonely edge & Other boundary edges & Rail edges\\
\midrule
$u_2$ & $x_2^{r+1}x_q^{r+1}$ & $x_2^{r+2}x_q^{r+2}$
& $\mathcal B_{0,r}$\\[1ex]
$x_3^0$ & $x_2^{r+1}x_q^{r+1}$ & $u_1u_2,\ x_2^{r+2}x_q^{r+2}$
& $\mathcal B_{0,r}$\\[1ex]
$x_t^{r+2}$ & $x_2^{r+1}x_q^{r+1}$ & $u_1u_3,\ u_2x_2^0$
& $\mathcal B_{0,r}$\\[1ex]
$u_1$ & $x_2^0x_3^0$ & $x_2^{r+1}x_2^{r+2},\ x_q^{r+2}x_t^{r+2}$
& $\mathcal A_{0,r}$\\[1ex]
$x_q^{r+1}$ & $x_2^0x_3^0$ & $u_2u_3,\ x_2^{r+2}x_t^{r+2}$
& $\mathcal A_{0,r}$\\[1ex]
$x_2^{r+2}$ & $x_2^0x_3^0$ & $u_2u_3$
& $\mathcal A_{0,r}$\\[1ex]
$\bar s_k$ & $x_2^0x_3^0$ & $u_2u_3,\ x_q^{r+1}x_q^{r+2},\ x_2^{r+1}x_2^{r+2}$
& $\mathcal A_{0,k-1}\cup\mathcal B_{k,r}$\\
\bottomrule
\end{tabular}
\end{center}

The degree verification is the same as in the preceding table:
alternating selections cover each root-free internal cell once,
and at cell $k$ the change from $\mathcal A$ to $\mathcal B$
gives both rail edges at $\bar s_k$, whose rung is supplied by
$\delta(\bar s_k)$. The boundary edges complete the required degrees.
Thus every row gives a join containing the stated nonincident lonely
edge. The rows cover all $r+6$ vertices outside $A$, so $A=\Act(G)$.

Finally, the two unique perfect matchings can be written as
\[
\begin{aligned}
M_{x_2^0x_3^0}
&=\{x_2^0x_3^0,u_2u_3,x_2^{r+1}x_2^{r+2},x_1^{r+2}x_3^{r+2}\}
  \cup\mathcal A_{0,r},\\
M_{x_2^{r+1}x_q^{r+1}}
&=\{u_1u_3,u_2x_2^0,x_2^{r+1}x_q^{r+1},x_2^{r+2}x_q^{r+2}\}
  \cup\mathcal B_{0,r}.
\end{aligned}
\]
At each $s_k$, these matchings use the two different rail edges and
neither uses the rung $s_k\bar s_k$.
At $u_3$, they use $u_2u_3$ and $u_1u_3$, respectively;
at $x_q^{r+2}$, they use its two different triangle edges.
Thus every indicated tip has the required matching pattern.
All arrow targets lie outside $A=\Act(G)$.
Together with the claim and the center path
$x_3^0-x_2^0-x_2^{r+1}-x_q^{r+1}$, this verifies all conditions
of \cref{def:unequal-standard-mark}.
\end{proof}

\begin{theorem}\label{thm:all-seeds-standard}
Every primitive $\U_2$ kernel admits at least one standard marking.
\end{theorem}

\begin{proof}
The ten finite primitive kernels are covered by
Lemma~\ref{lem:finite-seed-certificate}, and the four parity cases of
the two infinite kernel families are covered by
Lemmas~\ref{lem:infinite-common-mark}
and~\ref{lem:infinite-unequal-mark}.
\end{proof}

\subsection{Complete construction theorem}\label{sec:main-proof}

Only now do we deduce the concrete structure of an arbitrary graph in $\U_2$ from the facts that primitive kernels admit standard markings and that standard markings are closed under the three operations. In particular, the fact that the lonely edges of an equal-matching graph lie in triangles is proved here rather than assumed in Section~\ref{sec:three-operations}.

\begin{theorem}\label{thm:standard-mark-completeness}
The following hold.
\begin{enumerate}[label=\textup{(\roman*)}]
\item Every $G\in\U_2^{=}$ admits at least one equal-matching standard marking.
\item Every $G\in\U_2^{\ne}$ admits at least one unequal-matching standard marking.
\end{enumerate}
\end{theorem}

\begin{comment}
    More strongly, let
\[
C=G_0,G_1,\cdots,G_s=G,
\qquad
G_{i+1}=G_i^{v_i},\qquad \Sigma_{G_i}(v_i)=\varnothing
\]
be any zero-loss expansion sequence from a primitive kernel to $G$. If one chooses on $G_0$ the corresponding standard marking established in the previous section, then the marking can be propagated step by step along the whole sequence; every step is necessarily and exclusively one of Operations E, C, T compatible with the current marking.
\end{comment}

\begin{proof}
Take any $G\in\U_2$. By \cref{prop:kernel-reduction}, there is a primitive kernel $C$ and a zero-loss expansion sequence
\[
C=G_0,G_1,\cdots,G_s=G,
\qquad
G_{i+1}=G_i^{v_i},\qquad v_i\in\Act(G_i).
\]
By Lemma~\ref{lem:matching-type-preserved}, the matching type is preserved at every zero-loss expansion. Thus if $G\in\U_2^{=}$, then every $G_i$ and the primitive kernel $C$ are of equal-matching type; if $G\in\U_2^{\ne}$, the entire sequence is of unequal-matching type. By Theorem~\ref{thm:all-seeds-standard}, $G_0=C$ has a standard marking of the corresponding type. 

We now induct on $i$, assuming that $G_i$ carries the marking propagated along the sequence.

If $G_i$ is of equal-matching type, the definition of the standard marking gives
\[
\Act(G_i)=\{\text{the four endpoints of the two lonely edges}\}.
\]
Since $v_i\in\Act(G_i)$, the current step can only be Operation E at an endpoint of a lonely edge. Lemma~\ref{lem:operation-E} then gives an equal-matching standard marking on $G_{i+1}$.

If $G_i$ is of unequal-matching type, the definition gives
\[
\Act(G_i)=\{\text{the two center vertices}\}
\mathbin{\dot\cup}
\{\text{all tips}\}.
\]
If $v_i$ is a center vertex, the step is necessarily Operation C, and Lemma~\ref{lem:operation-C} supplies the standard marking on $G_{i+1}$. If $v_i$ is a tip, the step is necessarily Operation T, and Lemma~\ref{lem:operation-T} gives the next marking. These two cases exhaust all zero-loss vertices.

Induction to $i=s$ yields a standard marking on $G$ and simultaneously proves that every step of the given zero-loss sequence is identified exactly by one of the three operations.
\end{proof}

Let $\mathcal C_2=\mathcal E_2\cup\{D_0,D_1,D_2\}$ be the set of
primitive $\U_2$ kernels, and put
\[
\mathcal C_2^{=}=\mathcal C_2\cap\U_2^{=},\qquad
\mathcal C_2^{\ne}=\mathcal C_2\cap\U_2^{\ne}.
\]

\begin{theorem}\label{thm:main}
For every Klee graph $G$, the following statements hold.
\begin{enumerate}[label=\textup{(\roman*)}]
\item $G\in\U_2^{=}$ if and only if $G$ can be obtained from some $C\in\mathcal C_2^{=}$ by finitely many applications of endpoint advance at a lonely-edge endpoint;
\item $G\in\U_2^{\ne}$ if and only if $G$ can be obtained from some $C\in\mathcal C_2^{\ne}$ by finitely many applications of center expansion and tip extension;
\item 
\(
\U_2=\U_2^{=}\mathbin{\dot\cup}\U_2^{\ne}
\)
is exactly the union of the two classes produced by these constructions.
\end{enumerate}
At every step, the three vertices of the new triangle, the two lonely edges, the center vertices, and the tip directions are updated uniquely according to the labeling rules in Section~\ref{sec:three-operations}.
\end{theorem}

\begin{proof}
First, every graph constructed by the three operations belongs to the corresponding class of $\U_2$. Every seed kernel belongs to $\U_2$ and, by Theorem~\ref{thm:all-seeds-standard}, carries a standard marking. The root of each of Operations E, C, T belongs to the current zero-loss vertex set, so every step satisfies $\Sigma_G(v)=\varnothing$; by \cref{prop:loss-formula}, the expanded graph still has exactly two lonely edges. \cref{lem:operation-E}, \cref{lem:operation-C} and \cref{lem:operation-T} guarantee that the standard marking is updated according to the definitions. The matching type is preserved by Lemma~\ref{lem:matching-type-preserved}, so starting from $\mathcal C_2^{=}$ produces only graphs in $\U_2^{=}$, while starting from $\mathcal C_2^{\ne}$ produces only graphs in $\U_2^{\ne}$.

Conversely, Theorem~\ref{thm:standard-mark-completeness} shows that the reverse of any kernel-reduction sequence for any $G\in\U_2$ is realized step by step by Operations E, C, T. If $G\in\U_2^{=}$, only Operation E occurs throughout the sequence; if $G\in\U_2^{\ne}$, only Operations C and T occur. Hence all three assertions of Theorem~\ref{thm:main} follow.
\end{proof}

Part~(i) recovers the class characterized recursively
in~\cite[Theorem~4.12]{NMGK2026}, using the complete primitive-kernel
list obtained here. Parts~(ii) and~(iii) incorporate the unequal-matching
type into the construction.

\begin{corollary}\label{cor:equal-triangle-structure}
If $G\in\U_2^{=}$, then each of its two lonely edges lies in a triangle; the complete set of zero-loss expansion vertices is exactly the four endpoints of the two lonely edges, and those four vertices form a noninterfering set.
\end{corollary}

\begin{proof}
By Theorem~\ref{thm:standard-mark-completeness}, $G$ admits an equal-matching standard marking. The conclusion is precisely the content of Definition~\ref{def:equal-standard-mark} together with Lemma~\ref{lem:equal-mark-noninterference}.
\end{proof}

\begin{corollary}\label{cor:unequal-center-tip-structure}
If $G\in\U_2^{\ne}$, then it admits at least one marking of the following form: the two lonely edges together with one center edge form a three-edge center path; the complete set of zero-loss expansion vertices consists exactly of the two center vertices and a collection of directed tips; and this zero-loss vertex set is noninterfering. A center expansion creates one new tip, while a tip extension merely moves one tip.
\end{corollary}

\begin{proof}
Choose an unequal-matching standard marking from Theorem~\ref{thm:standard-mark-completeness}, and apply Definition~\ref{def:unequal-standard-mark}, \cref{lem:operation-C} and  \cref{lem:operation-T}.
\end{proof}

\begin{comment}
 \begin{corollary}\label{cor:explicit-procedure}
To generate graphs in $\U_2$ from any primitive kernel, it is enough to retain the following local marking data.
\begin{enumerate}
\item Equal-matching type: record the two lonely edges and the triangles containing them; one may expand any endpoint of a lonely edge and update according to Figure~\ref{fig:operation-E}.
\item Unequal-matching type: record the three-edge center path and the tip arrows; one may expand any center vertex or any arrow tail, updating according to Figures~\ref{fig:operation-C} and~\ref{fig:operation-T}, respectively.
\end{enumerate}
  New vertices are always labeled $v_x,v_y,v_z$ by their old neighbors. Section~\ref{sec:kernel-marking} provides the initial markings for every kernel, and the three closure lemmas govern arbitrarily many subsequent steps. Thus the construction is unambiguous at the label level.  
\end{corollary}   

\begin{remark}\label{rem:mark-nonunique}
Completeness of standard markings is an existence statement, not a uniqueness statement. The same unlabeled graph in $\U_2$ may have different primitive-kernel contraction endpoints and different zero-loss expansion sequences, and may therefore acquire different but valid center paths or tip markings. We prove only that every given zero-loss sequence can be propagated by the three operations; we do not claim uniqueness of the marking or the decomposition.
\end{remark}
\end{comment}

\section*{Conclusion}
In this paper, we completed the classification of 3-connected cubic graphs with at least two lonely edges by giving a full structural characterization of \(\mathcal{U}_2\), the family of 3-connected cubic graphs with exactly two lonely edges. We identified its minimally non-contractible ``primitive kernels'' and developed constructive operations that generate every graph in \(\mathcal{U}_2\) from these kernels. The equal-matching case recovers the recursive characterization of Narayana et al.~\cite{NMGK2026}; the center and tip operations cover the unequal-matching case. Together with the earlier characterization of \(\mathcal{U}_k\) for \(k\ge 3\), this yields a complete description of the classes \(\mathcal{U}_k\) for all \(k\ge 2\). A natural remaining question is how to characterize the structure of \(\mathcal{U}_1\).

\appendix

\clearpage
\end{document}